\documentclass[13pt,a4paper,leqno]{amsart}
\usepackage{xspace}

\usepackage{amsfonts,amssymb}
\usepackage{amsthm}
\usepackage[all]{xy}

\newcommand{\be}{\begin{otherlanguage}{english}}
\newcommand{\ee}{\end{otherlanguage}}

\theoremstyle{definition}
\newtheorem{defn}[subsection]{Definition}
\newtheorem{rem}[subsection]{Remark}
\theoremstyle{plain}
\newtheorem{lemma}[subsection]{Lemma}
\newtheorem{prop}[subsection]{Proposition}
\newtheorem{thm}[subsection]{Theorem}
\newtheorem*{thm*}{Th\'eor\`eme}
\newtheorem{cor}[subsection]{Corollary}

\theoremstyle{remark}

\numberwithin{equation}{subsection}

\newcommand{\beq}{\begin{equation}}
\newcommand{\eeq}{\end{equation}}

\newcommand{\ra}{\rightarrow}
\newcommand{\lra}{\longrightarrow}

\newcommand{\xra}{\xrightarrow}

\newcommand{\Z}{\mathbb{Z}}
\newcommand{\F}{\mathbb{F}}
\newcommand{\G}{\mathbb{G}}
\newcommand{\Q}{\mathbb{Q}}
\newcommand{\A}{\mathbb{A}}
\newcommand{\bP}{\mathbb {P}}

\newcommand{\bD}{\mathbf {D}}
\newcommand{\bK}{\mathbb {K}}

\newcommand{\ff}{\mathfrak {F}}
\newcommand{\fJ}{\mathfrak {J}}
\newcommand{\fK}{\mathcal {S}}
\newcommand{\fS}{\mathfrak {S}}
\newcommand{\fV}{\mathfrak {V}}

\newcommand{\cF}{\mathcal {F}}
\newcommand{\cG}{\mathcal {G}}
\newcommand{\calH}{\mathcal {H}}

\newcommand{\cO}{\mathcal {O}}
\newcommand{\cP}{\mathcal {P}}
\newcommand{\cQ}{\mathcal {Q}}

\newcommand{\rF}{\mathrm {F}}
\newcommand{\rU}{\mathrm {U}}
\newcommand{\rH}{\mathrm {H}}

\newcommand{\lam}{\lambda}

\newcommand{\OK}{\cO_K}
\newcommand{\Ab}{\overline{A}}
\newcommand{\Bb}{\overline{B}}
\newcommand{\Kb}{\overline {K}}
\newcommand{\Sb}{\overline{S}}
\newcommand{\fb}{\overline{f}}
\newcommand{\Gb}{\overline{G}}
\newcommand{\Hb}{\overline{H}}
\newcommand{\Xb}{\overline{X}}
\newcommand{\xb}{\overline{x}}
\newcommand{\OKb}{\cO_{\overline {K}}}
\newcommand{\m}{\mathfrak{m}}
\newcommand{\esb}{{\overline{s}}}
\newcommand{\et}{\mathrm{\acute{e}t}}

\newcommand{\ib}{{\overline{i}}}
\newcommand{\jb}{{\overline{j}}}
\newcommand{\kb}{{\overline{k}}}
\newcommand{\etab}{{\overline{\eta}}}
\newcommand{\afp}{\mathrm{AFP}_{\OK}}
\newcommand{\fppf}{\mathrm{fppf}}

\newcommand{\Ap}{{}_pA}
\newcommand{\Bp}{{}_pB}

\newcommand {\gk}{\cG_K}
\newcommand {\gns}{\gk{\text -}\mathrm{Ens}}
\newcommand {\Grs}{\mathrm{Gr}_S}
\newcommand {\resp}{\emph{resp.\xspace}}

\newcommand{\ie}{\emph{i.e.} }

\DeclareMathOperator{\Gal}{\mathrm{Gal}}
\DeclareMathOperator{\Coker}{\mathrm{Coker}} 
\DeclareMathOperator{\Ker}{\mathrm{Ker}}

\DeclareMathOperator{\im}{\mathrm{Im}}
\DeclareMathOperator{\Pic}{\mathrm{Pic}}

\DeclareMathOperator{\Lie}{\mathrm{Lie}}
\DeclareMathOperator{\Hom}{\mathrm Hom}

\DeclareMathOperator{\Ext}{\mathrm Ext}  

\DeclareMathOperator{\cHom}{\mathcal{H}\mathit{om}} 
\DeclareMathOperator{\cExt}{\mathcal{E}\mathit{xt}} 

\DeclareMathOperator{\Spec}{\mathrm {Spec}}

\begin{document}

\title{Canonical subgroups  of  Barsotti-Tate  groups }
\author{Yichao Tian}
\address{LAGA, Institut Galil\'ee, Universit\'e Paris 13,
93430 Villetaneuse, France}
\email{tian@math.univ-paris13.fr}
\date{}
\maketitle

\begin{abstract}Let $S$ be the spectrum of  a complete discrete valuation ring with fraction field
of characteristic $0$ and perfect residue field  of characteristic
 $p\geq 3$. Let $G$ be a truncated
 Barsotti-Tate group of level 1 over $S$.  If ``$G$ is not too supersingular'', a condition that
will be explicitly expressed in terms of the valuation of a certain
determinant,  we prove that we can
canonically lift the kernel
of the Frobenius endomorphism of its special fibre to a subgroup scheme of $G$,
finite and flat over $S$.
We call it  the canonical subgroup of $G$.
\end{abstract}

\section{Introduction}

\subsection{}
Let $\OK$ be a complete discrete valuation ring with fraction field $K$
of characteristic $0$ and perfect residue field $k$ of characteristic $p>0$.
We put $S=\Spec(\OK)$ and denote by $s$ (resp. $\eta$) its closed (resp.
generic) point. Let $G$ be  a truncated  Barsotti-Tate group  of
level $1$ over $S$.
If $G_s$ is ordinary, the kernel of its Frobenius endomorphism
is a multiplicative group scheme and can be uniquely lifted to a closed subgroup
scheme of $G$, finite and flat over $S$.
If we do not assume $G_s$ ordinary but
only that ``$G$ is not too supersingular'', a condition that
will be explicitly expressed in terms of the valuation of a certain
determinant, we will prove that we can still
canonically lift the kernel
of the Frobenius endomorphism of $G_s$ to a subgroup scheme of $G$,
finite and flat over $S$.
We call it {\em the canonical subgroup} of $G$.
Equivalently, under the same condition,
we will prove that the Frobenius endomorphism of $G_s$
can be canonically lifted to an isogeny of truncated Barsotti-Tate groups over $S$.
This problem was first raised by Lubin in 1967 and solved by himself for
$1$-parameter formal groups \cite{lubin}.
A slightly weaker question was asked by Dwork in 1969
for abelian schemes and answered also by him for elliptic curves \cite{dwork}:
namely, could we extend the construction of the
canonical subgroup in the ordinary case to a ``tubular neighborhood''
(without requiring that it lifts the kernel of the Frobenius)?
The dimension one case played a fundamental role
in the pioneering work of Katz on $p$-adic modular forms \cite{Kz}.
For higher dimensional abelian schemes, Dwork's conjecture
was first solved by Abbes and Mokrane \cite{AM}; our approach is
a generalization of their results.
Later, there have been other proofs, always for abelian
schemes, by Andreatta and Gasbarri \cite{AG}, Kisin and Lai
\cite{KL} and Conrad \cite{Co}.

\subsection{} For an $S$-scheme $X$,
we denote by $X_1$ its reduction modulo $p$. The valuation $v_p$ of
$K$, normalized by $v_p(p)=1$, induces a truncated valuation
$\cO_{S_1}\backslash \{0\}\ra \Q\cap[0,1)$. Let $G$ be a truncated
Barsotti-Tate group  of level $1$ and height $h$ over $S$, $G^\vee$
be its Cartier dual, and $d$ be the dimension of the Lie algebra of
$G_s$ over $k$. The Lie algebra $\Lie(G_1^\vee)$ of $G_1^\vee$ is a
free $\cO_{S_1}$-module of rank $d^*=h-d$, canonically isomorphic to
$\Hom_{(S_1)_{\fppf}}(G_1,\G_{a})$ (\cite{Il} 2.1). The Frobenius
homomorphism of $\G_a$ over $S_1$ induces an endomorphism $F$ of
$\Lie(G_1^\vee)$, which is semi-linear with respect to the Frobenius
homomorphism of $\cO_{S_1}$. We define the Hodge height
\eqref{hodge-height} of $\Lie(G_1^\vee)$ to be the truncated
valuation of the determinant of a matrix of $F$. This invariant
measures the ordinarity of $G$.

\subsection{} Following \cite{AM}, we construct the canonical subgroup
of a truncated Barsotti-Tate group over $S$
by the ramification theory of Abbes and Saito \cite{AS}. Let
$G$ be a commutative finite and flat group scheme over $S$.
In \cite{AM}, the authors defined a canonical exhaustive
decreasing filtration $(G^a,a\in\Q_{\geq0})$ by
finite, flat and closed subgroup schemes of $G$. For a real number $a\geq 0$,
we put $G^{a+}=\cup_{b>a}G^b$, where $b$ runs over rational numbers.

\begin{thm}\label{main-thm} Assume that $p\geq 3$,
and let $e$ be the absolute ramification index of $K$ and
$j=e/(p-1)$. Let $G$ be a truncated  Barsotti-Tate group  of level
$1$ over $S$, $d$ be the dimension of the Lie algebra
of $G_s$ over $k$. Assume that the Hodge height of
$\Lie(G_1^\vee)$ is strictly smaller than $1/p$. Then,

\emph{(i)} the subgroup scheme $G^{j+}$ of $G$ is locally free of rank
$p^d$ over $S$;

\emph{(ii)} the special fiber of $G^{j+}$ is the kernel
of the Frobenius endomorphism of $G_{s}$.
\end{thm}

\subsection{}
Statement (i) was proved by Abbes and Mokrane \cite{AM} for the
kernel of multiplication by $p$ of an abelian scheme over $S$. We
extend their result to truncated Barsotti-Tate groups by using a
theorem of Raynaud to embed $G$ into an abelian scheme over $S$. To
prove statement (ii), which we call {\em ``the lifting property of
the canonical subgroup''}, we give a new description of the
canonical filtration of a finite, flat and commutative group scheme
over $S$ killed by $p$ in terms of {\em congruence groups}. Let $\Kb$ be an algebraic closure of the fraction field of $S$, $\cO_{\Kb}$ be the integral closure of $\cO_K$ in $\Kb$. Put $\Sb=\Spec(\cO_{\Kb})$. For
every $\lam\in\OKb$ with $0\leq v_p(\lam)\leq 1/(p-1)$, Sekiguchi,
Oort and Suwa \cite{SOS} introduced a finite and flat group scheme
$G_\lam$ of order $p$ over $\Sb$ (see \eqref{defn-cong}); following
Raynaud, we call it the congruence group of level $\lambda$. If
$v_p(\lam)=0$,  $G_{\lam}$ is isomorphic to the multiplicative group
scheme $\mu_p=\Spec (\cO_{\Kb}[X]/\bigl(X^p-1)\bigr)$ over $\Sb$; and if $v_p(\lam)=1/(p-1)$, $G_{\lam}$ is
isomorphic to the  the constant \'etale group scheme $ \F_p$. For
general $\lam\in \cO_{\Kb}$ with $0\leq \lam\leq 1/(p-1)$, there is a
canonical $\Sb$-homomorphism $\theta_\lam:G_\lam\ra \mu_p$, such that $\theta\otimes \Kb$ is an isomorphism. For a finite, flat and commutative group scheme $G$
over $S$ killed by $p$, $\theta_\lam$ induces a homomorphism
\[
\theta_\lam(G)\colon \Hom_{\Sb}(G,G_{\lam}) \rightarrow
G^\vee(\Kb)=\Hom_{\Sb}(G,\mu_p).
\] We prove that it is injective, and its
image depends only on the valuation $a=v_p(\lam)$; we denote it by
$G^\vee(\Kb)^{[ea]}$, where $e$ is the absolute ramification index
of $K$ (the multiplication by $e$ will be justified later).
Moreover, we get a
 decreasing exhaustive filtration $(G^\vee(\Kb)^{[a]}, a\in \Q\cap[0,\frac{e}{p-1}])$.

\begin{thm}\label{main-thm2}
Let $G$ be a finite, flat and commutative group scheme over $S$ killed by $p$.
Under the canonical perfect pairing
\[
G(\Kb)\times G^\vee(\Kb)\ra \mu_p(\Kb),
\]
we have for any rational number $a\in \Q_{\geq 0}$,
\[
G^{a+}(\Kb)^{\perp}=
\begin{cases}G^\vee(\Kb)^{[\frac{e}{p-1}-\frac{a}{p}]},&
{\rm if}\  \ 0\leq a\leq \frac{ep}{p-1},\\
G^\vee(\Kb),&{\rm if} \ \ a> \frac{ep}{p-1}.\end{cases}\]
\end{thm}
Andreatta and Gasbarri \cite{AG} have used congruence
groups to prove the existence of the canonical subgroup for
abelian schemes. This theorem explains the relation between the
approach via the ramification theory of \cite{AM} and this paper,  and
the approach of \cite{AG}.

\subsection{} This article is organized as follows.
For the convenience of the reader, we recall in section 2 the theory
of ramification of group schemes over a complete discrete valuation
ring, developed in \cite{AS} and \cite{AM}. Section 3 is a summary
of the results in \cite{AM} on the canonical subgroup of an abelian
scheme over $S$. Section 4 consists of  some preliminary results on
the $\fppf$ cohomology of  abelian schemes. In  section 5, we define
the Bloch-Kato filtration for a  finite, flat and commutative group
scheme over $S$ killed by $p$. Using this filtration, we prove
Theorem \ref{main-thm}(i) in section 6. Section 7 is dedicated to
the proof of Theorem \ref{main-thm2}. Finally in section 8, we
complete the proof of Theorem \ref{main-thm}(ii).

\subsection{} This article is a part of the author's
thesis at Universit\'e Paris 13. The author would like to express
his great gratitude to his thesis advisor Professor A. Abbes for
leading him to this problem and for his helpful comments on earlier
versions of this work. The author thanks Professors W. Messing and
M. Raynaud for their  help. He is also grateful to the referee for
his careful reading and very valuable comments.

\subsection{Notation}\label{notations}
In this article, $\OK$ denotes a complete discrete valuation ring
with fraction field $K$ of characteristic $0$, and residue field $k$
of characteristic $p>0$. Except in Section \textbf{2}, we will
assume that $k$ is perfect. Let $\Kb$ be an algebraic closure of
$K$, $\gk=\Gal(\Kb/K)$ be the Galois group of $\Kb$ over $K$, $\OKb$
be the integral closure of $\OK$ in $\Kb$, $\m_{\Kb}$ the maximal
ideal of $\cO_{\Kb}$, and  $ \overline{k}$ be the residue field of
$\cO_{\Kb}$.

We put $S=\Spec(\OK)$, $\Sb=\Spec(\OKb)$, and denote by $s$ and
$\eta$ (\resp $\esb$ and $\etab$) the closed and generic point of
$S$ (\resp of $\Sb$) respectively.

We fix a uniformizer $\pi$ of $\OK$. We will use two valuations $v$
and $v_p$ on $\OK$, normalized respectively by $v(\pi)=1$ and
$v_p(p)=1$; so we have $v=e v_p$, where $e$ is the absolute
ramification index of $K$. The valuations $v$ and $v_p$ extend
uniquely to $\Kb$; we denote the extensions also by $v$ and $v_p$.
For a rational number $a\geq 0$, we put $\m_a=\{x\in \Kb;v_p(x)\geq
a\}$ and $\Sb_a=\Spec (\cO_{\Kb}/\m_a)$. If $X$ is a scheme over
$S$, we will denote respectively by $\Xb$, $X_{\esb}$ and $\Xb_a$
the schemes obtained by base change of $X$ to $\Sb$, $\esb$ and
$\Sb_a$.

If $G$ is a commutative, finite and flat group scheme over $S$, we
will denote by $G^\vee$ its  Cartier dual. For an abelian scheme $A$
over $S$, $A^\vee$ will denote the dual abelian scheme, and ${}_pA$
the kernel of multiplication by $p$, which is a finite and flat
group scheme over $S$.

\section{Ramification theory of finite flat group schemes over $S$}

\subsection{}
We begin by recalling the main construction of \cite{AS}. Let $A$ be
a finite and flat $\cO_K$-algebra. We fix a  finite presentation of
$A$ over $\OK$
\[
0\ra I\ra \OK[x_1,\cdots,x_n]\ra A\ra 0,
\]
or equivalently, an $S$-closed immersion of $i\colon
\Spec(A)\rightarrow \A^n_{S}$. For a rational number $a> 0$, let
$X^a$ be the tubular neighborhood of $i$ of thickening $a$
(\cite{AS} Section 3, \cite{AM} 2.1). It is an affinoid subdomain of
the $n$-dimensional closed unit disc over $K$ given by
\[
X^a(\Kb)=\{(x_1,\cdots,x_n)\in \OKb^n|\ \ v(f(x_1,\cdots,x_n))\geq
a, \ \ \ \ \forall f\in I\}.
\]
Let $\pi_0(X^a_{\Kb})$ be the set of geometric connected components
of $X^a$. It is a finite $\gk$-set that does not depend on the
choice of the presentation (\cite{AS} Lemma 3.1). We put
\begin{equation}\label{form-Fa}
\cF^a(A)=\pi_0(X_{\Kb}^a).
\end{equation}
For two rational numbers $b\geq a> 0$, $X^b$ is
an affinoid sub-domain of $X^a$. So there is a natural
transition map $\cF^b(A)\ra \cF^a(A)$.

\subsection{} We denote by $\afp$ the category of finite
flat $\OK$-algebras, and by $\gns$ the category of finite sets with
a continuous action  of $\gk$. Let
\begin{eqnarray*}
\cF\colon \afp^{\circ}&\ra &\gns\\
A&\mapsto& \Spec(A)(\Kb)
\end{eqnarray*}
be the functor of geometric
points. For $a\in \Q_{>0}$, \eqref{form-Fa} gives rise to a
functor
\[
\cF^a: \afp^\circ \ra \gns.
\]
For $b\geq a\geq 0$, we have morphisms of
functors $\phi^a:\cF\ra \cF^a$ and $\phi^a_b: \cF^b\ra \cF^a$, satisfying
 the relations $\phi^a=\phi^b\circ\phi^a_b$ and
$\phi^a_c=\phi^a_b\circ\phi^b_c$ for $c\geq b\geq a\geq 0$.
To stress the dependence on $K$, we will denote $\cF$ (\resp $\cF^a$)
by $\cF_K$ (\resp $\cF_K^a$).
These functors behave well only for finite, flat and
relative complete intersection algebras over $\OK$ (EGA IV 19.3.6).
We refer to \cite{AS} Propositions 6.2 and 6.4 for their main properties.

\begin{lemma}[\cite{AM} Lemme 2.1.5]\label{ext-inv} Let $K'/K$ be an extension
(not necessarily finite) of complete discrete valuation fields with
ramification index $e_{K'/K}$. Let $A$ be a finite, flat and
relative complete intersection algebra over $\OK$. Then
we have a canonical isomorphism $\cF^{ae_{K'/K}}_{K'}(A')\simeq
\cF^a_K(A)$ for all $a\in \Q_{>0}$.
\end{lemma}

\subsection{} Abbes and Saito show that the projective system of functors
$(\cF^a,\cF\ra \cF^a)_{a\in \Q_{\geq0}}$ gives rise
to an exhaustive decreasing filtration ($\gk^a, a\in \Q_{\geq0}$)
of the group $\gk$, called the \emph{ ramification  filtration}
(\cite{AS} Proposition 3.3). Concretely, if
$L$ is a finite Galois extension of $K$ contained in $\Kb$,
$\Gal(L/K)$ is the Galois group of $L/K$, then the quotient
filtration $(\Gal(L/K)^a)_{a\in \Q_{\geq0}}$ induced by
$(\gk^a)_{a\in\Q_{\geq 0}}$ is determined by the following canonical
isomorphisms
\[
\cF^a(L)\simeq \Gal(L/K)/\Gal(L/K)^a.
\]
For a real number $a\geq 0$, we put
$\gk^{a+}=\overline{\cup_{b>a}\gk^{b}}$, and if $a>0$
$\gk^{a-}=\cap_{b<a}\gk^b$ , where $b$ runs over rational numbers.
Then $\gk^{0+}$ is the inertia subgroup of $\gk$ (\cite{AS}
Proposition 3.7).

\subsection{} We recall the definition of the canonical filtration
of a finite and flat group schemes over $S$, following \cite{AM}.
Let $\Grs$ be the category of finite, flat and commutative group schemes  over $S$.
Let $G$ be an object of $\Grs$ and $a\in
\Q_{\geq 0}$.  Then there is a natural group structure on $\cF^a(A)$
(\cite{AM} 2.3), and the canonical surjection $\cF(A)\ra\cF^a(A)$ is
a homomorphism of $\gk$-groups. Hence, the kernel
$G^a(\Kb)=\Ker(\cF(A)\ra \cF^a(A))$ defines a subgroup scheme
$G^a_{\eta}$ of $G_\eta$ over $\eta$, and the schematic closure
$G^a$ of $G^a_{\eta}$ in $G$ is a closed subgroup scheme of $G$,
locally free of finite type over $S$. We put $G^0=G$. The exhaustive
decreasing filtration $(G^a,a\in \Q_{\geq0})$ defined above is
called \emph{the canonical filtration of $G$} (\cite{AM} 2.3.1).
Lemma \ref{ext-inv} gives immediately the following.

\begin{lemma}\label{gp-ext-inv}
Let $K'/K$ be an extension (not necessarily finite) of complete
discrete valuation fields with ramification index $e_{K'/K}$, $\cO_{K'}$ be
the ring of integers  of $K'$, and $\overline{K'}$ be an algebraic
closure of $K'$ containing $\Kb$. Let $G$ be an object of $\Grs$ and
$G'=G\times_S\Spec(\cO_{K'})$. Then we have a canonical isomorphism
$G^a(\Kb)\simeq G'^{ae_{K'/K}}(\overline{K'})$ for all $a\in \Q_{>0}$.
\end{lemma}

\subsection{} For an object $G$ of $\Grs$ and $a\in \Q_{\geq 0}$,
we denote $G^{a+}=\cup_{b>a}G^b$ and if $a>0$,
$G^{a-}=\cap_{0<b<a}G^b$, where $b$ runs over rational numbers. The
construction of the canonical filtration is functorial: a morphism
$u:G\ra H$ of $\Grs$ induces canonical homomorphisms $u^a: G^a\ra
H^a$, $u^{a+}:G^{a+}\ra H^{a+}$ and  $u^{a-}:G^{a-}\ra H^{a-}$.

\begin{prop}[\cite{AM} Lemmes 2.3.2 and  2.3.5]\label{comp-con}
\emph{(i)} For any object $G$ of $\Grs$, $G^{0+}$ is the neutral
connected component of $G$.

\emph{(ii)} Let $u: G\ra H$ be a finite flat and surjective morphism in
$\Grs$ and $a\in \Q_{>0}$. Then the homomorphism $u^a(\Kb):
G^a(\Kb)\ra H^a(\Kb)$ is surjective.
\end{prop}

\subsection{} Let $A$ and $B$ be two abelian schemes over $S$, $\phi:A\ra
B$ be an isogeny (\ie a finite flat morphism of group schemes), and
$G$ be the kernel of $\phi$. Let $\nu$ (\resp $\mu$) be the generic
point of the special fiber $A_s$ (\resp $B_s$), and
$\widehat{\cO}_{\nu}$ (\resp $\widehat{\cO}_{\mu}$) be the
completion of the local ring of $A$ at $\nu$ (\resp of $B$ at
$\mu$).  Let $M$ and $L$ be the fraction fields of
$\widehat{\cO}_{\nu}$ and $\widehat{\cO}_{\mu}$ respectively. So we
have the cartesian diagram
\[
\xymatrix{\Spec M\ar[r]\ar[d]&\Spec
\widehat{\cO}_{\nu}\ar[r]\ar[d]&A\ar[d]\\
\Spec L\ar[r]&\Spec \widehat{\cO}_{\mu}\ar[r]&B}
\]
We fix a separable closure $\overline{L}$ of $L$ containing $\Kb$,
and an imbedding of $M$ in $\overline{L}$. Since  $\phi:A\ra B$ is a
$G$-torsor, $M/L$ is a Galois extension, and we have a canonical
isomorphism \beq\label{gp-gal}
G(\Kb)=\cF_{K}(G)\stackrel{\sim}{\rightarrow}
\cF_L(\widehat{\cO}_{\nu})=\Gal(M/L). \eeq Using the same arguments
of (\cite{AM} 2.4.2), we prove the following

\begin{prop}\label{ram-fini-abel} For all rational numbers
$a\geq 0$, the isomorphism \emph{(\ref{gp-gal})}
induces an isomorphism $G^a(\Kb)\simeq \Gal(M/L)^a$.
\end{prop}

\section{Review of the abelian scheme case following \cite{AM}}

From this section on, we assume that the residue field $k$ of $\OK$
is perfect of characteristic  $p>0$.

\subsection{} Let $X$ be a smooth and proper scheme over $S$, and $\overline{X}=X\times_S\Sb$. We consider the
cartesian diagram
\[\xymatrix{X_{\esb}\ar[r]^{\ib}\ar[d]&\overline{X}\ar[d]&X_{\etab}\ar[d]\ar[l]_{\jb}\\
\esb=\Spec \kb\ar[r]&\Sb&\etab=\Spec \Kb\ar[l]}\] and the sheaves of
$p$-adic vanishing cycles on $X_{\esb}$
\beq\label{cycle-evan}\Psi^q_X=\ib^*R^q\jb_*(\Z/p\Z(q)),\eeq where
$q\geq 0$ is an integer and $\Z/p\Z(q)$ is the Tate twist of
$\Z/p\Z$. It is clear that $\Psi^0_X\simeq \Z/p\Z$. By the base change theorem for proper morphisms,
 we have a spectral sequence
\beq\label{suite-spec}E_2^{p,q}(X)=\rH^p(X_\esb,\Psi^q_X)(-q)\Longrightarrow
\rH^{p+q}(X_{\etab},\Z/p\Z),\eeq which induces an exact sequence
\beq\label{cycle-abel}0\ra \rH^1(X_{\esb},\Z/p\Z)\ra
\rH^1(X_{\etab},\Z/p\Z)\xra{u} \rH^0(X_\esb, \Psi^1_X)(-1)\ra
\rH^2(X_\esb, \Z/p\Z).\eeq

\subsection{} The Kummer's exact sequence $0\ra \mu_p\ra \G_m\ra
\G_m\ra 0$ on $X_{\etab}$
 induces the symbol map
\beq\label{cycle-symbol}h_{\overline{X}}:\ib^*\jb_*\cO_{X_{\etab}}^\times\ra\Psi^1_X.\eeq
We put $\rU^0\Psi^1_X=\Psi^1_X$, and for $a\in \Q_{>0}$,
\beq\label{BK-cycle}\rU^a\Psi^1_X=h_{\overline{X}}(1+\pi^{a}\ib^*\cO_{\overline{X}}),\eeq
where by abuse of notation $\pi^a$ is  an element in $\cO_{\Kb}$ with
$v(\pi^a)=a$ .  We have $\rU^{a}\Psi^1_X=0$ if $a\geq
\frac{ep}{p-1}$ (\cite{AM} Lemme 3.1.1).

Passing to the cohomology, we get a filtration on
$\rH^1(X_{\etab},\Z/p\Z)$ defined by:
\begin{align}\rU^0\rH^1(X_{\etab},\Z/p\Z)&=\rH^1(X_{\etab},\Z/p\Z),\nonumber\\
\rU^a\rH^1(X_{\etab},\Z/p\Z)&=u^{-1}(\rH^0(X_\esb,\rU^a\Psi^1_X)(-1)),\quad
\text{for $a\in \Q_{>0}$,}\label{BK-abel}
\end{align}
called the \emph{ Bloch-Kato filtration}.

\begin{thm}[\cite{AM} Th\'eor\`eme 3.1.2]\label{Thm-BK-abel} Let $A$ be an abelian scheme over $S$,
$\Ap$  its kernel of multiplication by $p$, and $e'=\frac{ep}{p-1}$.
Then under the canonical perfect pairing
\begin{equation}\label{pairing-3.3}\Ap(\Kb)\times \rH^1(A_{\etab},\Z/p\Z)\ra \Z/p\Z,\end{equation} we
have for all $a\in \Q_{\geq0}$,
\[\Ap^{a+}(\Kb)^{\perp}=\begin{cases}\rU^{e'-a}\rH^1(A_\esb,\Z/p\Z)\quad &\text{if $0\leq a\leq e';$}\\
\rH^1(A_{\etab},\Z/p\Z)\quad &\text{if $a>e'$.}\end{cases}\]
\end{thm}

\subsection{}
 Let $X$ be a scheme over
$S$, $\overline{X}=X\times_S \Sb$. For all $a\in \Q_{>0}$, we put
$\Sb_a=\Spec (\OKb/\m_a)$ and $\Xb_a=X\times_S \Sb_a$
\eqref{notations}. We denote by $\bD((\Xb_1)_{\et})$ the derived
category of abelian \'etale sheaves over $\Xb_1$.  A morphism of
schemes is called \emph{syntomic}, if it is flat and of complete
intersection.

 Let $X$ be a syntomic  and quasi-projective $S$-scheme, $r$
and $n$ be  integers with $r\geq 0$ and $n\geq 1$. In \cite{Ka},
Kato constructed
 a canonical object $\fJ^{[r]}_{n,\Xb}$ in
$\bD((\Xb_1)_{\et})$,  and if $0\leq r\leq p-1$ a morphism
$\varphi_r:\fJ^{[r]}_{n,\Xb}\ra \fJ^{[0]}_{n,\Xb}$, which can be
roughly seen as ``$1/p^r$ times of the Frobenius map''. We refer to
\cite{Ka} and (\cite{AM} 4.1.6) for details of these constructions.
Let $\fK_n(r)_{\Xb}$ be the fiber cone of the morphism $\varphi_r-1$
for $0\leq r\leq p-1$; so we have a distinguished triangle in
$\bD((\Xb_1)_{\et})$ \beq\label{syn-triangle}\fK_n(r)_{\Xb}\ra
\fJ^{[r]}_{n,\Xb}\xra{\varphi_r-1} \fJ^{[0]}_{n,\Xb}\xra{+1}.\eeq
The complexes $\fK_n(r)_{\Xb}$ $(0\leq r\leq p-1)$ are called   the
\emph{ syntomic complexes} of $\Xb$. For our purpose, we recall here
some of their properties for $r=1$.

\subsection{} According to (\cite{Ka} section I.3),  for any integer $n\geq 1$, there exists a
surjective symbol map
\[ \cO_{\overline{X}_{n+1}}^\times\ra \calH^1(\fK_n(1)_{\Xb}).\]
For a geometric point $\xb$ of $X_{\esb}$, we put
$$\fS^1_{\xb}=\biggl(\cO_{\Xb,\xb}[\frac{1}{p}]\biggr)^\times=(\ib^*\jb_*\cO_{X_{\etab}}^\times)_{\xb}.$$
By (\cite{Ka} I.4.2), the above symbol map at $\xb$  factorizes
through the canonical surjection $\cO_{\Xb,\xb}^{\times}\ra
\fS_{\xb}^1/p^n\fS^1_{\xb}$.

\begin{thm}[\cite{Ka} I.4.3]\label{isom-syn-cycle} Assume that $p\geq 3$. Let $X$ be a
smooth and quasi-projective scheme over $S$.  Then there is a
canonical isomorphism $\calH^1(\fK_1(1)_{\Xb})\tilde{\ra}\Psi^1_X$,
which is compatible with the symbol maps $\fS^1_{\xb}\ra
\calH^1(\fK_1(1)_{\Xb})_{\xb}$ and $h_{\Xb}:\fS^1_{\xb}\ra
\Psi^1_{X,\xb}$ \emph{(\ref{cycle-symbol})}.
\end{thm}

\subsection{} Let $X$  be a smooth and quasi-projective scheme over
$S$.  Let $\phi_{\Xb_1}$ and $\phi_{\Sb_1}$ be the absolute
Frobenius morphisms of $\Xb_1$ and $\Sb_1$, and let $\Xb_1^{(p)}$ be
the scheme defined by the cartesian diagram
\[\xymatrix{\Xb_1^{(p)}\ar[d]\ar[r]^{w}&\Xb_1\ar[d]\\
\Sb_1\ar[r]^{\phi_{\Sb_1}}&\Sb_1.}\] For all integers $q\geq 0$, we
denote by $\rF$ the composed morphism
\[\Omega^q_{\Xb_1/\Sb_1}\xra{w^*}\Omega^q_{\Xb_1^{(p)}/\Sb_1}\ra \Omega^q_{\Xb_1/\Sb_1}/d(\Omega^{q-1}_{\Xb_1/\Sb_1}),\]
where the second morphism is induced by the Cartier isomorphism
\[C^{-1}_{\Xb_1/\Sb_1}:\Omega^q_{\Xb_1^{(p)}/\Sb_1}\xra{\sim}\calH^q(\Omega^{\bullet}_{\Xb_1/\Sb_1}).\]
Let $c$ be the class in $\cO_{\Sb_1}$ of a $p$-th root of $(-p)$. We
set
\beq\label{defn-cp}\cP=\Coker\biggl(\cO_{\Xb_1}\xra{\rF-c}\cO_{\Xb_1}\biggr),\eeq
\[\cQ=\Ker\biggl(\Omega^1_{\Xb_1/\Sb_1}\xra{\rF-1}\Omega^1_{\Xb_1/\Sb_1}/d(\cO_{\Xb_1})\biggr).\]

\begin{prop}[\cite{AM} 4.1.8]\label{syn-prop} The notations are those as above, and we assume moreover that $p\geq
3$. Let $\xb$ be a geometric point in $X_{\esb}$.

\emph{(i)} There exist canonical isomorphisms
\[\cP\xra{\sim}\Coker\biggl(\calH^0(\fJ^{[1]}_{1,\Xb})\xra{\varphi_1-1}\calH^0(\fJ^{[0]}_{1,\Xb})\biggr),\quad
\cQ\xra{\sim}
\Ker\biggl(\calH^1(\fJ^{[1]}_{1,\Xb})\xra{\varphi_1-1}\calH^1(\fJ^{[0]}_{1,\Xb})\biggr),\]
so  the distinguished triangle \emph{(\ref{syn-triangle})} gives
rise to an exact sequence \beq\label{syn-suite} 0\ra
\cP\xra{\alpha}\calH^1(\fK_1(1)_{\Xb})\xra{\beta}\cQ\ra0.\eeq

\emph{(ii)} Let $e(T)=\sum_{i=0}^{p-1}T^i/i!\in \Z_p[T]$.  Then the
morphism $\alpha_{\xb}$ in \eqref{syn-suite} is induced by the map
$\cO_{\Xb_1,\xb}\ra\fS_{\xb}^1/p\fS^1_{\xb}$ given by
\[a\mapsto e(-\tilde{a}(\zeta-1)^{p-1}),\] where $\tilde{a}$ is a
lift of $a\in \cO_{\Xb_1,\xb}$ and $\zeta\in \Kb$ is a primitive
$p$-th root of  unity.

\emph{(iii)} The composed map
\[\fS_{\xb}^1/p\fS^1_{\xb}\xra{\mathrm{symbol}}\calH^1(\fK_1(1)_{\Xb})_{\xb}\xra{\beta}
\cQ_{\xb}\] is the unique morphism sending $a\in
\cO_{\Xb,\xb}^\times$ to $a^{-1}da\in \cQ_{\xb}.$
\end{prop}

\begin{rem}\label{rem-syn-BK} Statement (ii) of Proposition \ref{syn-prop} implies that,
via the canonical isomorphism $\calH^1(\fK_1(1)_{\Xb})\simeq
\Psi^1_X$ (\ref{isom-syn-cycle}),
 $\cP$ can be identified with the  submodule $\rU^e\Psi^1_X$ of $\Psi^1_X$ defined in (\ref{BK-cycle}).
\end{rem}

\begin{prop}[\cite{AM} 4.1.9]\label{prop-fc} Assume that $p\geq 3$. Let $X$ be a smooth
projective scheme over $S$, $t=(p-1)/p$.

\emph{(i)} The morphism $\rF-c:\cO_{\Xb_1}\ra \cO_{\Xb_1}$
factorizes through the quotient morphism $\cO_{\Xb_1}\ra
\cO_{\Xb_t}$, and we have an exact sequence \beq\label{fc-suite}
0\ra \F_p\ra \cO_{\Xb_t}\xra{\rF-c}\cO_{\Xb_1}\ra \cP\ra 0. \eeq

\emph{(ii)} Let $\delta_E: \rH^0(\Xb_1,\cP)\ra \rH^2(\Xb_1,\F_p)$ be
the cup-product with the class of $(\ref{fc-suite})$ in
$\Ext^2(\cP,\F_p)$, and
\[d_2^{0,1}:\rH^0(\Xb_\esb,\Psi^1_X)(-1)\ra \rH^2(\Xb_{\esb},\F_p)\]
be the connecting morphism  in \emph{(\ref{cycle-abel})}.  Then the
composed morphism
\[\rH^0(\Xb_1,\cP)\ra \rH^0(\Xb_{\esb},\calH^1(\fK_1(1)_{\Xb}))\xra{\sim}\rH^0(\Xb_{\esb},\Psi^1_X)\xra{d_2^{0,1}(1)}\rH^2(\Xb_{\esb},\F_p)(1)\]
coincides with  $\zeta\delta_E$, where  the middle isomorphism is
given by Theorem \ref{isom-syn-cycle}, and $\zeta $ is a chosen
$p$-th root of  unity.

\emph{(iii)} Assume moreover that
$\rH^0(\Xb_r,\cO_{\Xb_r})=\cO_{S_r}$ for  $r=1$ and $t$. Then we
have an exact sequence \beq\label{fc-coh-suite} 0\ra
\rH^1(X_{\esb},\F_p)\ra\Ker\biggl(\rH^1(\Xb_1,\cO_{\Xb_t})\xra{\rF-c}\rH^1(\Xb_1,\cO_{\Xb_1})\biggr)\ra
\rH^0(\Xb_1,\cP)\xra{\delta_E}\rH^2(X_{\esb},\F_p) \eeq

\end{prop}

\subsection{}\label{hodge-height}

Let $M$ be  a free $\cO_{\Sb_1}$-module of rank $r$ and $\varphi:
M\ra M$ be a semi-linear endomorphism  with respect to the absolute
Frobenius of $\cO_{\Sb_1}$. Following \cite{AM}, we call $M$ a
\emph{$\varphi$-$\cO_{\Sb_1}$-module} of rank $r$. Then $\varphi(M)$
is an $\cO_{\Sb_1}$-submodule of $M$, and there exist rational
numbers $0\leq a_1\leq a_2\leq \cdots\leq a_r\leq 1$, such that
\[M/\varphi(M)\simeq \oplus_{i=1}^{r}\cO_{\Kb}/\m_{a_i}.\]
We define the \emph{Hodge height} of $M$ to be $\sum_{i=1}^ra_i$.
For any rational number $0\leq t\leq 1$, we put
$M_t=M\otimes_{\cO_{\Sb_1}}\cO_{\Sb_t}.$

\begin{prop}[\cite{AG} 9.1  and  9.7]\label{prop-HW} Assume that $p\geq  3$.
Let $\lam$ be an element in $\cO_{\Kb}$, and $r\geq 1$  an
integer. We assume $v=v_p(\lam)<\frac{1}{2}$ and  let $M$ be a
$\varphi$-$\cO_{\Sb_1}$-module of rank $r$ such that its Hodge
height is strictly smaller than $v$.

\emph{(i)} The morphism $\varphi-\lam:M\ra M$ factorizes through the
canonical map $M\ra M_{1-v}$ and the kernel  of $
\varphi-\lam:M_{1-v}\ra M$ is an $\F_p$-vector space of dimension
$r$.

\emph{(ii)} Let $N_0$ be the kernel of the morphism $M_{1-v}\ra M$
induced by $\varphi-\lam$, and $N$ be the  $\cO_{\Kb}$-submodule of
$M_{1-v}$ generated by $N_0$. Then we have
$\dim_{\overline{k}}(N/\m_{\Kb}N) =\dim_{\F_p}N_0=r$.

\end{prop}



\subsection{}  We can now summarize the strategy of \cite{AM} as follows.
 Let $A$ be a projective abelian  scheme of dimension $g$ over $S$. By
Proposition \ref{prop-fc}(iii),  we have
\begin{align}\label{equ-fc-abel}\dim_{\F_p}\rH^1(\Ab_1,\F_p)&+\dim_{\F_p}\Ker
\biggl(\rH^0(\Ab_1,\cP)\xra{\delta_E} \rH^2(\Ab_1,\F_p)\biggr)\\
&=\dim_{\F_p}\Ker\biggl(\rH^1(\Ab_1,\cO_{\Ab_t})\xra{\rF-c}
\rH^1(\Ab_1,\cO_{\Ab_1})\biggr).\nonumber
\end{align}
By Remark \ref{rem-syn-BK} and Proposition \ref{prop-fc}(ii), the
left hand side equals $\dim_{\F_p}\rU^e\rH^1(A_{\etab},\Z/p\Z)$
(\ref{BK-abel}). Taking account of  Theorem \ref{Thm-BK-abel}, we
get
\beq\label{dim-can-abel}2g-\dim_{\F_p}\bigl(\Ap^{j+}(\Kb)\bigr)=\dim_{\F_p}\Ker\biggl(\rH^1(\Ab_1,\cO_{\Xb_t})\xra{\rF-c}
\rH^1(\Ab_1,\cO_{\Xb_1})\biggr),\eeq where $j=\frac{e}{p-1}$.
Applying  \ref{prop-HW}(i)  to $M=\rH^1(\Ab_1,\cO_{\Xb_1})$ and
$\lam=c$, we obtain immediately the first statement of Theorem
\ref{main-thm} for projective abelian schemes. In fact, Abbes and
Mokrane  proved a less optimal bound on the Hodge height (\cite{AM}
5.1.1).

\section{Cohomological preliminaries }

\subsection{} Let $f:X\ra T$ be a  proper,
flat and finitely presented morphism of schemes.
We work with the $\fppf$-topology on $T$, and denote by $\Pic_{X/T}$
the relative Picard functor $R^1_{\fppf}f_*(\G_m)$.

If $T$ is the spectrum of a field,
 $\Pic_{X/T}$ is representable by a group scheme locally  of finite type over $k$.
We denote by $\Pic^0_{X/T}$ the neutral component,
and put $\Pic^\tau_{X/T}=\bigcup_{n\geq 1}n^{-1}\Pic^0_{X/T}$,
where $n:\Pic_{X/T}\ra \Pic_{X/T}$ is the multiplication by $n$. Then $\Pic^0_{X/T}$
 and $\Pic^\tau_{X/T}$ are open sub-group schemes of $\Pic_{X/T}$.

For a general base,
$\Pic_{X/T}$ is representable by an algebraic space over $T$ (\cite{Ar} thm. 7.3).
  We denote by $\Pic^0_{X/T}$ (\resp $\Pic^{\tau}_{X/T}$) the subfunctor of $\Pic_{X/T}$
which consists of all elements whose restriction to all fibres $X_t$, $t\in T$,
 belong to $\Pic^0_{X_t/t}$ (\resp $\Pic^{\tau}_{X_t/t}$). By (SGA 6 XIII, thm 4.7), the canonical inclusion
  $\Pic^{\tau}_{X/T} \ra \Pic_{X/T}$ is relatively
representable by an open  immersion.

\subsection{} Let  $f: A\ra T$ be an abelian scheme.
 If $T$ is the  spectrum of a field, the N\'eron-S\'everi group $\Pic_{A/T}/\Pic^0_{A/T}$ is torsion free, \ie we have $\Pic^0_{A/T}=\Pic^{\tau}_{A/T}$.  This coincidence remains true for a general base $T$ by the definitions of $\Pic^0_{A/T}$ and $\Pic^\tau_{A/T}$.
  Moreover,
  $\Pic^\tau_{A/T}$ is formally smooth (cf. \cite{Mu} Prop. 6.7),
  and $\Pic^\tau_{A/T}$ is actually open and closed in $\Pic_{A/T}$, and is representable by a
  proper and smooth algebraic space over $T$, \ie an abelian
  algebraic space over $T$. By a theorem of
   Raynaud (\cite{FC} Ch. 1, thm. 1.9), every abelian algebraic
  space over $T$ is automatically an abelian scheme over $T$. So
  $\Pic^0_{A/T}=\Pic^\tau_{A/T}$ is an abelian scheme,
 called the \emph{dual abelian scheme of $A$}, and denoted
  by $A^\vee$.

Let $H$ be a commutative group scheme locally free of finite type
over $T$.  Recall the following isomorphism due to Raynaud
(\cite{Ray} 6.2.1):
\beq\label{im-abel}R^1_{\fppf}f_*(H_A)\xra{\sim}\cHom(H^\vee,{\Pic}_{A/T}),
\eeq where $H_A=H\times_T A $, $H^\vee$ is the Cartier dual  of $H$,
and $\cHom$ is taken for  the $\fppf$-topology on $T$.

\begin{prop}\label{assump4.1} Let $A$ be an abelian scheme over a scheme $T$, and $H$
a commutative group scheme locally free of finite type over $T$.
Then we have  canonical isomorphisms
\begin{align}
\cExt^1(A,H)&\xra{\sim}\cHom(H^\vee,A^\vee),\label{ext-abel} \\
\Ext^1(A,H)&\xra{\sim}\rH^0_{\fppf}(T,\cExt^1(A,H))\xra{\sim}\Hom(H^\vee,A^\vee).\label{ext-abel-global}
\end{align}
\end{prop}

\begin{proof}For any $\fppf$-sheaf $E$ on $T$, we have
$$\cHom(H^\vee, \cHom(E,\G_m))\simeq \cHom(H^\vee\otimes_{\Z} E,\G_m)\simeq
\cHom(E,H).$$ Deriving this isomorphism of functors in $E$ and
putting $E=A$, we obtain a spectral sequence
$$E_2^{p,q}=\cExt^p(H^\vee, \cExt^q(A,\G_m))\Longrightarrow\cExt^{p+q}(A,H).$$
Since  $\cHom(A,\G_m)=0$, the exact sequence
$$0\ra E_2^{1,0}\ra \cExt^1(A,H)\ra E_2^{0,1}\ra E_2^{2,0}$$
induces an isomorphism $$\cExt^1(A,H)\simeq
\cHom(H^\vee,\cExt^1(A,\G_m)).$$ Then \eqref{ext-abel} follows from
the canonical identification $\cExt^1(A,\G_m)\simeq A^\vee$ (\cite{Dem}  2.4). For
 \eqref{ext-abel-global}, the spectral sequence
\[E^{p,q}_2=\rH^p_{\fppf}(T,\cExt^q(A,H))\Longrightarrow\Ext^{p+q}(A,H)\]
induces a long exact sequence
\[0\ra \rH^1_{\fppf}(T,\cHom(A,H))\ra \Ext^1(A,H)\xra{(1)} \rH^0_{\fppf}(T,\cExt^1(A,H))\ra \rH^2_{\fppf}(T,\cHom(A,H)).\]
Since $\cHom(A,H)=0$, the  arrow $(1)$ is an isomorphism, and
\eqref{ext-abel-global} follows by applying the functor
$\rH^0_{\fppf}(T,\_)$ to \eqref{ext-abel}.
\end{proof}

\subsection{} The  assumptions are those of \eqref{assump4.1}. We define a canonical morphism
\beq\label{morph-ext-tor}\Ext^1(A,H)\ra \rH^1_{\fppf}(A,H)\eeq as
follows. Let $a$ be an element in $ \Ext^1(A,H)$ represented by the
extension $0\ra H\ra E\ra A\ra 0$. Then the $\fppf$-sheaf $E$ is
representable by a scheme over $T$, and is naturally a $H$-torsor
over $A$. The image of $a$ by the homomorphism \eqref{morph-ext-tor}
is defined to be the class of the torsor $E$. Since this
construction is functorial in $T$,  by passing to sheaves, we obtain
a canonical morphism \beq\label{varphi-H} \cExt^1(A,H)\ra
R_{\fppf}^1f_{*}(H_A).\eeq Via the isomorphisms (\ref{im-abel}) and
(\ref{ext-abel}), we check that \eqref{varphi-H} is induced by the
canonical map $A^\vee=\Pic^0_{A/T}\ra \Pic_{A/T}$.

Since $H$ is faithfully flat and finite over $T$, the inverse image
of the $\fppf$-sheaf $H$ by $f$ is representable by $H_A$, \ie we
have $f^*(H)=H_A$. Therefore, we deduce an adjunction morphism
\begin{equation}\label{morph-adjoint-finite}H\ra R^0_{\fppf}f_*(H_A).\end{equation}

\begin{prop}\label{cor-im-ext}  Let  $f: A\ra T$ be an abelian scheme,
and $H$ be a commutative group scheme locally free of finite type
over $T$. Then the canonical maps \eqref{varphi-H} and
\eqref{morph-adjoint-finite} are isomorphisms.
\end{prop}

\begin{proof} First, we prove that \eqref{varphi-H} is an isomorphism. By (\ref{im-abel}) and
(\ref{ext-abel}), we have to verify that the canonical morphism
$$\cHom(H^\vee, A^\vee)\ra \cHom(H^\vee, \Pic_{A/T})$$ is an
isomorphism.  Let $g: H^\vee\ra \Pic_{A/T}$ be a homomorphism over
$T$. For every  $t\in T$, the induced morphism $g_t: H^\vee_t\ra
\Pic_{A_t/t}$
  falls actually in $\Pic^\tau_{A_t/t}$, because $H^\vee$ is a finite group scheme.
Hence, by the definition of $\Pic^{\tau}_{A/T}$,
 the homomorphism $g$ factorizes through the canonical inclusion
$A^\vee=\Pic^\tau_{A/T}\ra \Pic_{A/T}$; so the canonical morphism
$\Hom(H^\vee, A^\vee)\ra \Hom(H^\vee, \Pic_{A/T})$ is an
isomorphsim.

Secondly, we prove that \eqref{morph-adjoint-finite} is an
isomorphism. For  $T$-schemes $U$ and  $G$, we denote $G_U=G\times_T
U$. We must verify that for any $T$-scheme $U$, the adjunction
morphism
\begin{equation}\label{fppf-iso}
\varphi(U):H(U)\ra R_{\fppf}^0f_*(H_A)(U)=H(A_U)\end{equation}
 is an isomorphism. We note that
$H(U)=H_U(U)$ and $H(A_U)=H_U(A_U)$; therefore, up to taking base
changes, it suffices to prove that  $\varphi(T)$ \eqref{fppf-iso} is
an isomorphism. We remark that $f$ is surjective, hence $\varphi(T)$
is injective. To prove the surjectivity of $\varphi(T)$, we take an
element $h\in H(A)$, \ie a morphism of $T$-schemes $h:A\ra H$; by
rigidity lemma for abelian schemes (cf. \cite{Mu} Prop. 6.1), there
exists a section $s:T\ra H$ of the structure morphism $H\ra T$ such
that $s\circ f=h$. Hence we have $\varphi(T)(s)=h$, and $\varphi(T)$
is an isomorphism.
\end{proof}



\begin{cor}\label{rem-im-ext} Let $T$ be the spectrum of a stirctly henselian local ring,
$f:A\ra T$  an abelian scheme, and $H$  a finite \'etale group
scheme over $T$. Then we have canonical isomorphisms
\[\rH^1_{\et}(A,H)\simeq\rH^1_{\fppf}(A,H)\simeq \Ext^1(A,H).\]

\end{cor}
\begin{proof}The first isomorphism follows from the \'etaleness of $H$
(\cite{Gr}, 11.7). For the second one, the ``local-global'' spectral
sequence induces a long exact sequence
\[0\ra \rH^1_{\fppf}(T,R_{\fppf}^0f_*(H_A))\ra \rH^1_{\fppf}(A,H)\ra \rH^0_{\fppf}(T,R_{\fppf}^1f_*(H_A))
\ra \rH^2_{\fppf}(T,R^0_{\fppf}f_*(H_A)).\] By Prop.
\ref{cor-im-ext}, we have $R_{\fppf}^0f_*(H_A)=H$. Since  $T$ is
strictly henselian and $H$ \'etale, we have
$\rH^q_{\fppf}(T,H)=\rH^q_{\et}(T,H)=0$ for all integers $q\geq 1$.
Therefore, we obtain $\rH^1_{\fppf}(A,H)\xra{\sim}
\rH^0_{\fppf}(T,R_{\fppf}^1f_*(H_A))$, and the corollary follows
from  \ref{cor-im-ext} and \eqref{ext-abel-global}.
\end{proof}
\subsection{} Let $T$ be a scheme, and $G$ be a commutative  group scheme locally free of finite type over $T$.
We denote by $\G_a$  the additive group scheme, and  by $\Lie
(G^\vee)$ the Lie algebra of $G^\vee$. By Grothendieck's duality
formula (\cite{MM} II \S 14), we have a canonical isomorphism
\beq\label{Lie-alg-fini}\Lie (G^\vee)\simeq \cHom_{T}(G,\G_a),\eeq
where we have regarded $G$ and $\G_a$ as abelian  $\fppf$-sheaves on
$T$.
 If $T$ is of characteristic $p$ and
$G$ is a truncated Barsotti-Tate group over $T$, then $\Lie(G^\vee)$ is a
locally free of finite type $\cO_T$-module (\cite{Il} 2.2.1(c)).

 Similarly, for an abelian scheme $f:A\ra T$, we
have a canonical isomorphism (\cite{BBM} 2.5.8)
\beq\label{Lie-alg-abel}\Lie(A^\vee)\simeq \cExt^1_{T}(A,\G_a)\simeq
R^1f_{*}(\G_a).\eeq  In the sequel, we will frequently use the
identifications \eqref{Lie-alg-fini} and \eqref{Lie-alg-abel}
without any indications.\\

The following Lemma is indicated by W. Messing.

\begin{lemma}\label{phi-module} Let $L$ be an algebraically closed field of characteristic $p>0$,
 $R$ be an $L$-algebra integral over $L$,  and $M$ be
 a  module of finite presentation  over $R$, equipped with an
endomorphism $\varphi$ semi-linear with respect to the Frobenius of
$R$. Then the map $\varphi-1:M\ra M$ is surjective.
\end{lemma}

\begin{proof} First, we reduce the lemma to the case $R=L$.
 Consider $R$ as a filtrant inductive limit of finite $L$-algebras
$(R_i)_{i\in I}$.  Since $M$ is of finite presentation,  there
exists an $i\in I$, and an $R_i$-module $M_i$ of finite presentation
endowed with a Frobenius semi-linear endomorphism $\varphi_i$, such
that $M=M_i\otimes_{R_i}R$ and $\varphi=\varphi_i\otimes \sigma$,
where $\sigma$ is the Frobenius on $R$. For $j\geq i$, we put
$M_j=M_i\otimes_{R_i}R_j$ and
$\varphi_j=\varphi_i\otimes_{R_i}\sigma_j$, where $\sigma_j$ is the
Frobenius of $R_j$. In order to prove $\varphi-1$ is surjective on
$M$, it is sufficient to prove the surjectivity of $\varphi_j-1$ on each $M_j$ for $j\geq i$.
Therefore, we may assume that
$R$ is a finite dimensional $L$-algebra, and  $M$ is thus a finite
dimensional vector space over $L$.

We put $M_1=\bigcup_{n\geq 1}\Ker (\varphi^n)$ and
$M_2=\bigcap_{n\geq 1}\im (\varphi^n)$. Then we have a decomposition
$M=M_1\oplus M_2$ as $\varphi$-modules, such that $\varphi$ is
nilpotent on $M_1$ and bijective on $M_2$ (Bourbaki, Alg\`ebre VIII
\S 2 $n^{\circ}$ 2 Lemme 2). Therefore, it is sufficient to prove
the surjectivity of $\varphi-1$ in the following two cases:

(i) $\varphi$ is nilpotent. In this case, the endomorphism
$1-\varphi$ admits an inverse $1+\sum_{n\geq 1}\varphi^n$. Hence it
is surjective.

(ii) $\varphi$ is invertible.   We choose a basis of $M$ over $L$,
and let $U=(a_{i,j})_{1\leq i,j \leq n}$ be the matrix of $\varphi$
in this basis. The problem reduces to prove that the
 equation system
 $\sum_{j=1}^na_{i,j} x_j^p-x_i=b_i$ $(1\leq i \leq n)$ in $X=(x_1,\cdots,x_n)$
 has solutions for all $b=(b_1,\cdots,b_n)\in L^{n}$.
Since  $U$ is invertible, let $V=(c_{i,j})_{1\leq i,j\leq n}$ be its
inverse. Then the equation system $\sum_{j=1}^na_{i,j}
x_j^p-x_i=b_i$ is equivalent to $x_i^p-\sum_{j=1}^n c_{i,j}x_j=b'_i$
for $1\leq i\leq n$ with $b'=\sum_j c_{i,j} b_j$. But these $n$
equations define a finite \'etale cover of $\Spec L$ of degree
$p^n$. Hence they admit solutions in $L$, since $L$ is separably
closed. This completes the proof.
\end{proof}

\begin{cor}\label{ext2nul} Let $H$ be a Barsotti-Tate group or an abelian scheme over $\Sb_1$ \eqref{notations}. Then $\Ext^2(H,\F_p)=0$ for  the $\fppf$ topology on $\Sb_1$.
\end{cor}
\begin{proof}
 Let $K_0$ be the fraction field of the ring of Witt vectors with coefficients in $k$;
so $K$ is a finite extension of degree $e$ of $K_0$. Let
$\cO_{K_0}^{ur}$ be the ring of integers of the maximal unramified
extension of $K_0$ in $\Kb$. Then $\cO_{\Sb_1}=\OKb/p\OKb$ is
integral over the algebraically closed field
$\overline{k}=\cO^{ur}_{K_0}/p\cO^{ur}_{K_0}$. As $\Ext^2(H,\G_a)=0$
(\cite{BBM} Proposition 3.3.2), the Artin-Schreier's exact sequence
$0\ra\F_p\ra \G_a\xra{\rF-1}\G_a\ra0$  induces an exact sequence
\[\Ext^1(H,\G_a)\xra{\varphi-1}\Ext^1(H,\G_a)\ra \Ext^2(H,\F_p)\ra0.\]
 Since $\Ext^1(H,\G_a)$ is a free $\cO_{\Sb_1}$-module
(\cite{BBM} 3.3.2.1), the corollary follows immediately from Lemma
\ref{phi-module}.
\end{proof}

\section{The Bloch-Kato filtration for finite flat group schemes killed by $p$ }

\subsection{}\label{thm-raynaud} Recall the following theorem  of  Raynaud (\cite{BBM} 3.1.1):
\emph{Let $T$ be a scheme, $G$ be a
commutative group scheme locally free of finite type  over $T$. Then locally for
the Zariski topology, there exists a projective abelian scheme $A$ over $T$,
such that $G$ can be identified to  a closed subgroup of $A$.}

In particular, if $G$ is a commutative finite and flat group scheme
over $S=\Spec (\OK)$, we have an exact sequence of abelian
$\fppf$-sheaves over $S$ \beq\label{imm-abel}0\ra G\ra A\ra B\ra
0,\eeq where $A$ and $B$ are  projective abelian schemes over $S$.
In the sequel, such an exact sequence is called a \emph{resolution
of $G$ by abelian schemes}.

\subsection{}\label{functoriality} Let $f:X\ra Y$ be a morphism of proper and smooth $S$-schemes.
For any integer $q\geq 0$, we have  a base change morphism
$f_{\esb}^*(\Psi^q_Y)\ra \Psi^q_X$
 of $p$-adic vanishing cycles \eqref{cycle-evan} (SGA 7 XIII 1.3.7.1). For $q=1$, this morphism  respects the Bloch-Kato
filtrations
 (\ref{BK-cycle}), that is, it sends
 $f_{\esb}^*(\rU^a\Psi^1_Y)$ to $\rU^a\Psi^1_X$ for all $a\in \Q_{\geq 0}$.

 Passing to cohomology, we get
 a functorial map $
\rH^p(Y_{\esb}, \Psi^q_Y)(-q)\ra \rH^p(X_{\esb},\Psi^q_X)(-q)$ for
each pair of integers $p,q\geq 0$. These morphisms piece together to give a morphism
of spectral sequences \eqref{suite-spec} $E^{(p,q)}_2(Y)\ra
E^{(p,q)}_2(X)$, which converges to the  map
$\rH^{p+q}(Y_{\etab},\Z/p\Z)\ra \rH^{p+q}(X_{\etab},\Z/p\Z)$ induced
by $f_{\etab}^*$.
 Therefore, we have the following commutative diagram
\beq\label{suite-cycles1}\xymatrix{0\ar[r]&\rH^1(Y_{\esb},\Z/p\Z)\ar[r]\ar[d]^{\alpha_1}&\rH^1(Y_{\etab},\Z/p\Z)
\ar[r]\ar[d]^{\alpha_2}&\rH^0(Y_{\esb},\Psi^1_Y)(-1)\ar[d]^{\alpha_3}\ar[r]^{\quad \, d_2^{1,0}}&\rH^2(Y_\esb,\Z/p\Z)\ar[d]^{\alpha_4}\\
0\ar[r]&\rH^1(X_{\esb},\Z/p\Z)\ar[r]&\rH^1(X_{\etab},\Z/p\Z)\ar[r]&\rH^0(X_{\esb},\Psi^1_X)(-1)\ar[r]^{\quad
\, d_2^{1,0}}&\rH^2(X_{\esb},\Z/p\Z).} \eeq It is clear that the
Bloch-Kato filtration on $\rH^1(X_{\etab},\Z/p\Z)$ \eqref{BK-abel}
is functorial in $X$. More precisely, the following diagram is
commutative: \beq\label{funct-BK}\xymatrix{\rU^a\rH^1(Y_\etab,
\Z/p\Z)\ar@{^(->}[r]\ar[d]&\rH^1(Y_\etab,
\Z/p\Z)\ar[d]^{\alpha_2}\\
\rU^a\rH^1(X_\etab,\Z/p\Z)\ar@{^(->}[r]&\rH^1(X_\etab,\Z/p\Z)}\eeq

\subsection{}\label{subsect5.3} Let $G$ be a commutative finite and flat  group scheme over
 $S$ killed by $p$, and $0\ra G\ra A\ra B\ra 0$  a resolution  of
 $G$ by abelian schemes \eqref{imm-abel}. We apply the  construction \eqref{suite-cycles1} to the morphism $ A\ra B$.
 Using Corollary \ref{rem-im-ext}, we obtain immediately
that
\begin{align*}
\Ker \alpha_2&=\Ker \biggl(\Ext^1(B_\etab,\F_p)\ra \Ext^1(A_{\etab},\F_p)\biggr)\\
&=\Hom(G_{\etab},\F_p)=G^\vee(\Kb)(-1),\\
\Ker \alpha_1&=\Ker\biggl(\Ext^1(B_{\esb},\F_p)\ra \Ext^1(A_{\esb},\F_p)\biggr)\\
&=\Hom(G_{\esb},\F_p)=(G_{\et})^\vee(\Kb)(-1),
\end{align*}
where $G_{\et}=G/G^{0+}$  is the \'etale part of $G$ (cf.
\ref{comp-con}). Setting $N=\Ker \alpha_3$, we can complete
\eqref{suite-cycles1} as follows:
\beq\label{suite-cycles}\xymatrix{0\ar[r]
&(G_{\et})^\vee(\Kb)(-1)\ar[r]\ar@{^(->}[d]^{\gamma_1}&G^\vee(\Kb)(-1)\ar[r]^{u}\ar@{^(->}[d]^{\gamma_2} &N\ar@{^(->}[d]^{\gamma_3}\ar@^{-->}[r]&0\ar@^{-->}[d]\\
0\ar[r]&\rH^1(B_{\esb},\Z/p\Z)\ar[r]\ar[d]^{\alpha_1}&\rH^1(B_{\etab},\Z/p\Z)
\ar[r]\ar[d]^{\alpha_2}&\rH^0(B_{\esb},\Psi^1_B)(-1)\ar[d]^{\alpha_3}\ar[r]^{\quad  d_2^{1,0}(B)}&\rH^2(B_\esb,\Z/p\Z)\ar[d]^{\alpha_4}\\
0\ar[r]&\rH^1(A_{\esb},\Z/p\Z)\ar[r]&\rH^1(A_{\etab},\Z/p\Z)\ar[r]&\rH^0(A_{\esb},\Psi^1_A)(-1)\ar[r]^{\quad
 d_2^{1,0}(A)}&\rH^2(A_{\esb},\Z/p\Z).}\eeq We will show later
that the morphism $u$ is surjective.

\begin{defn}\label{filt-BK-fini} The assumptions are those of \ref{subsect5.3}.
 We call the \emph{Bloch-Kato filtration}  on
$G^\vee(\Kb)$, and denote  by $(\rU^aG^\vee(\Kb),a\in\Q_{\geq0})$,
the
 decreasing and exhaustive filtration defined by
$\rU^0G^\vee(\Kb)=G^\vee(\Kb)$,
 and   for $a\in \Q_{>0},$
\begin{equation}\label{BK-fini-abel}\rU^aG^\vee(\Kb)=\gamma_2^{-1}(\rU^a\rH^1(B_{\etab},\Z/p\Z))(1).\end{equation}

\end{defn}

\begin{prop}\label{prop-BK} Let $e'=\frac{ep}{p-1}$, $G$ be a commutative finite and flat group scheme
 over $S$ killed by $p$, and $0\ra G\ra A\ra B\ra 0$ be a resolution of $G$ by abelian schemes \eqref{imm-abel}.

\emph{(i)} For all $a\in \Q_{\geq 0}$, we have
\begin{align}\label{4.9.1}\rU^aG^\vee(\Kb)&\simeq\Ker \biggl(\rU^a\rH^1(B_\etab,
\Z/p\Z)(1)\xra{\alpha_2(1)}
\rU^a\rH^1(A_\etab,\Z/p\Z)(1)\biggr)\\
&\simeq u^{-1}\biggl(N(1)\cap
\rH^0(B_{\esb},\rU^{a}\Psi^1(B))\biggr),\nonumber\end{align} where
$N(1)$ is identified to a subgroup of $\rH^0(B_{\esb},\Psi^1(B))$ by
$\gamma_3(1)$ in \eqref{suite-cycles}.

 \emph{(ii)}  The morphism $u:
G^\vee(\Kb)(-1)\ra N$ in \eqref{suite-cycles} is surjective. In
particular, $N$ is contained in the kernel of the morphism
$d_2^{1,0}(B)$ in \eqref{suite-cycles}.

 \emph{(iii)} Under the
canonical perfect pairing \beq\label{pairing}G(\Kb)\times
G^\vee(\Kb)\ra \mu_p(\Kb),\eeq we have, for all $a\in \Q_{\geq 0},$
\[G^{a+}(\Kb)^{\perp}=\begin{cases}\rU^{e'-a}G^\vee(\Kb)\quad&\text{if $0\leq a\leq e'$;}\\
G^\vee(\Kb)\quad &\text{if $a>e'$.}\end{cases}\] In particular, the
filtration $(\rU^aG^\vee(\Kb),a\in\Q_{\geq 0})$ does not depend on
the resolution of $G$ by abelian schemes.
\end{prop}

\begin{proof} Statement (i) is obvious from  definition \ref{filt-BK-fini}
 and  diagrams \eqref{funct-BK}
and \eqref{suite-cycles}.

For (ii) and (iii), thanks to Lemma \ref{gp-ext-inv}, we need only
to prove the proposition after a base change $\cO_{K}\ra \cO_{K'}$,
where $K'/K$ is a finite extension. Therefore, up to such a base
change, we may add the following assumptions.

(1) We may assume that $k$ is algebraically closed,   $K$ contains a
primitive $p$-th root of unity, and $G(\Kb)=G(K)$.

(2) For $X=A$ or $B$, we consider the cartesian diagram
\[\xymatrix{X_{s}\ar[r]^{i}\ar[d]&X\ar[d]&X_{\eta}\ar[d]\ar[l]_j\\
s=\esb\ar[r]&S&\eta\ar[l]}\] and the \'etale sheaf
$\Psi^1_{X,K}=i^*R^1j_*(\Z/p\Z)$ over $X_s$. By an argument as in
the proof of (\cite{AM} 3.1.1), we may assume
 that $\rH^0(X_s,\Psi^1_X)=\rH^0(X_s,\Psi^1_{X,K})$.

 Since
$\rU^{a}\Psi^1_X=0$ for $a\geq e'$ (\cite{AM} Lemme 3.1.1), we have
$\rU^{e'}G^\vee(\Kb)=\Ker (u)(1)=(G_{\et})^\vee(\Kb).$ Statement
(iii) for  $a=0$ follows immediately from Proposition
\ref{comp-con}(i). The pairing (\ref{pairing}) induces a perfect
pairing \beq\label{ind-pair} G^{0+}(\Kb)\times \im(u)(1)\lra
\mu_p(\Kb).\eeq In particular, we have
$\dim_{\F_p}\bigl(\im(u)(1)\bigr)=\dim_{\F_p}\bigl(G^{0+}(\Kb)\bigr)$.

 Let $\mu$ (\resp $\nu$) be
 the generic point of $B_s=B_{\esb}$ (\resp of $A_s=A_{\esb}$), and $\overline{\nu}$
   be  a geometric point over $\nu$. Then $\overline{\nu}$ induces by the morphism $\nu\ra\mu$ a geometric
   point over $\mu$, denoted by $\overline{\mu}$.
   Let $\cO_\mu$ (\resp $\cO_\nu$) be the local ring of $B$ at $\mu$
(\resp of $A$ at $\nu$), $\cO_{\overline{\mu}}$ (\resp
$\cO_{\overline{\nu}}$) be the henselization of $\cO_{\mu}$ at
$\overline{\mu}$ (\resp of $\cO_{\nu}$ at $\overline{\nu}$). We
denote by $\widehat{\cO}_{{\mu}}$ and  $\widehat{\cO}_{{\nu}}$ the
completions of $\cO_{\mu}$ and $\cO_{\nu}$, and by $
(\widehat{\cO}_{\mu})^{\mathrm{h}}_{\overline{\mu}}$ (\resp
 by $(\widehat{\cO}_{\nu})^{\mathrm{h}}_{\overline{\nu}}$) the henselization  of $\widehat{\cO}_{\mu}$
(\resp of $\widehat{\cO}_{\nu}$) at $\overline{\mu}$ (\resp at
$\overline{\nu}$). Let $L_0$ (\resp $M_0$) be the fraction field of
$\cO_{\mu}$ (\resp of $\cO_{\nu}$), $L_0^{ur}$ (\resp $M_0^{ur}$) be
the fraction field of ${\cO}_{\overline{\mu}}$ (\resp of
$\cO_{\overline{\nu}}$). We denote by  $L$ (\resp by $M$)  the
fraction field of $\widehat{\cO}_{{\mu}}$ (\resp of
$\widehat{\cO}_{{\nu}}$),  and by $L^{ur}$ (\resp by $M^{ur}$) the
fraction field of
$(\widehat{\cO}_{\mu})^{\mathrm{h}}_{\overline{\mu}}$ (\resp of
$(\widehat{\cO}_{\nu})^{\mathrm{h}}_{\overline{\nu}}$). We notice
that $\pi$ is a uniformizing element in $\widehat{\cO}_{\mu}$ and in
$(\widehat{\cO}_{\mu})^{\mathrm{h}}_{\overline{\mu}}$.
\[\xymatrix{\cO_{\overline{\mu}}\ar[r]&L_0^{ur}\ar[r]&L^{ur}&(\widehat{\cO}_{\mu})^{\mathrm{h}}_{\overline{\mu}}\ar[l]\\
\cO_{\mu}\ar[u]\ar[r]&L_0\ar[r]\ar[u]&L\ar[u]&\widehat{\cO}_{\mu}\ar[l]\ar[u]}\]

Since  $G=\Ker(A\ra B)$, we have an identification $\Gal(M/L)\simeq
G(K)= G(\Kb)$. We fix a separable closure $\overline{L}$ of $L$, and
an imbedding of $M$ in $\overline{L}$, which induces a surjection
$\varphi: \Gal(\overline{L}/L)\ra \Gal(M/L)$.  By
\ref{ram-fini-abel}, we have
$\varphi(\Gal(\overline{L}/L)^a)=\Gal(M/L)^a=G^a(\Kb)$,  for all
$a\in \Q_{>0}$. In particular, we have
\beq\label{galois-surj}\varphi(\Gal(\overline{L}/L^{ur}))=\Gal(M/M\cap
L^{ur})=\Gal(M^{ur}/L^{ur})=G^{0+}(\Kb).\eeq Since $L^{ur}/L$ is
unramified, we have, for all $a\in \Q_{>0}$,
\beq\label{4.7.1}\Gal(M^{ur}/L^{ur})^{a}= \Gal(M/L)^{a}=G^{a}.\eeq

 Let $\rho_B$ be the composition of
the canonical morphisms
$$\rH^0(B_s,\Psi^1_B)=\rH^0(B_s,\Psi_{B,K}^1)\ra (\Psi_{B,K}^1)_{\overline{\mu}}=
\rH^1(\Spec L_0^{ur},\mu_p)\ra \rH^1(\Spec L^{ur},\mu_p),$$ and we
define $\rho_A$ similarly. By functoriality and \eqref{galois-surj},
we have a commutative diagram
\beq\label{inf-res}\xymatrix{0\ar[r]&N(1)\ar[r]\ar@{^(->}[d]&\rH^0(B_s,\Psi^1_B)\ar@{^(->}[d]^{\rho_B}\ar[r]&\rH^0(A_s,\Psi^1_A)\ar@{^(->}[d]^{\rho_A}\\
0\ar[r]&\rH^1(G^{0+}(\Kb),\mu_p)\ar[r]^{inf}&\rH^1(\Spec
L^{ur},\mu_p)\ar[r]^{res} &\rH^1(\Spec M^{ur},\mu_p),}\eeq where the
lower horizontal  row is the ``inflation-restriction'' exact
sequence in Galois cohomology.  By (\cite{BK} Prop. 6.1), $\rho_B$
and $\rho_A$ are injective. Hence  $\im(u)(1)\subset N(1)$ is
identified with a subgroup of $\rH^1(G^{0+}(\Kb),\mu_p)$. By
assumption (1), we have
$\rH^1(G^{0+}(\Kb),\mu_p)=\Hom(G^{0+}(\Kb),\mu_p(\Kb))$, which has
the same dimension over $\F_p$ as $G^{0+}(\Kb)$.  Hence by the
remark below \eqref{ind-pair}, we get
\beq\label{im-coincide}\im(u)(1)=N(1)\simeq
\rH^1(G^{0+}(\Kb),\mu_p)=\Hom(G^{0+}(\Kb),\mu_p(\Kb)).\eeq This
proves statement (ii) of the proposition.

Statement (iii) for  $a=0$ has been proved above. Since the
filtration $G^{a}$ is decreasing and $\rU^0G^\vee(\Kb)=G^\vee(\Kb)$,
we may assume, for the proof of (iii), that $0< a\leq e'$. It
suffices to prove
 that,  under the pairing \eqref{ind-pair}, we have
 \[G^{a+}(\Kb)^{\perp}=
\im(u)(1)\cap \rH^0(B_{s},\rU^{e'-a}\Psi^1_B).\]
 From \eqref{inf-res} and \eqref{im-coincide}, we check easily that
 \eqref{ind-pair} is identified with the canonical pairing
\[G^{0+}(\Kb)\times \rH^1(G^{0+}(\Kb),\mu_p)\ra \mu_p(\Kb).\]
Hence we are reduced to prove that, under this pairing, we have
\beq\label{ind-dual-fini}G^{a+}(\Kb)^{\perp}=
\rH^1(G^{0+}(\Kb),\mu_p)\cap \rH^0(B_{s},\rU^{e'-a}\Psi^1_B),\eeq
where the ``$\cap$'' is taken in $\rH^1(\Spec L^{ur},\mu_p)$
\eqref{inf-res}.

Let $h: (L^{ur})^\times\ra H^1(\Spec L^{ur},\mu_p)$ be the symbol
morphism. We define a decreasing filtration on $\rH^1(\Spec
L^{ur},\mu_p)$ in a  similar way  as (\ref{BK-cycle}) by
\[\rU^0\rH^1=\rH^1,
\quad \text{and}\quad
\rU^b\rH^1=h(1+\pi^b(\widehat{\cO}_{\mu})^{\mathrm{h}}_{\overline{\mu}})\text{
for all integers $b>0$.}\] We extend this definition to all
$b\in\Q_{\geq 0}$ by setting $\rU^{b}\rH^1=\rU^{[b]}\rH^1$, where
$[b]$ denotes the integer part of $b$. By (\cite{AM} Lemme 3.1.3),
for all $b\in\Q_{>0}$, we have
$\rH^0(B_s,\rU^b\Psi^1_B)=\rho_B^{-1}\biggl(\rU^{b}\rH^1(\Spec
L^{ur},\mu_p)\biggr).$ Therefore, the right hand side of
\eqref{ind-dual-fini} is
\beq\label{subgroup}\rH^1(G^{0+}(\Kb),\mu_p)\cap
\rU^{(e'-a)}\rH^1(\Spec L^{ur},\mu_p).\eeq
 We identify $\rH^1(G^{0+}(\Kb),\mu_p)$ with the character group
$\Hom(\Gal(M^{ur}/L^{ur}),\mu_p(\Kb))$;  then by (\cite{AM} Prop.
2.2.1), the subgroup \eqref{subgroup} consists of the characters
$\chi$ such that $\chi(\Gal(M^{ur}/L^{ur})^{(e'-a)+})=0.$ In view of
(\ref{4.7.1}), we obtain immediately (\ref{ind-dual-fini}), which
completes the proof.

\end{proof}

\begin{rem}The proof of \ref{prop-BK}(iii) follows the same strategy
as the proof of  (\cite{AM} Th\'eorem 3.1.2). The referee point out 
that we can also reduce \ref{prop-BK}(iii) to \ref{Thm-BK-abel}.
Indeed, the commutative diagram
\[\xymatrix{0\ar[r]&G\ar[r]\ar[d]^{\times p}&A\ar[r]\ar[d]^{\times p}&B\ar[r]\ar[d]^{\times p}&0\\
0\ar[r]&G\ar[r]&A\ar[r]&B\ar[r]&0}\] induces, by  snake lemma, an
exact sequence of finite group schemes $0\ra G\ra
\Ap\xra{\phi}\Bp\xra{\psi}G\ra0.$  We consider the following perfect
pairing
\[\xymatrix{&\Ap(\Kb)\ar[d]_{\phi}&\Ap(\Kb)\times \rH^1(A_{\etab},\F_p)\ar[rr]^-(0.5){(\bullet, \bullet)_A}&&\F_p\ar@{=}[d]&\rH^1(A_{\etab},\F_p)\\
\Bp^{a+}(\Kb)\ar@{^(->}[r]\ar@{->>}[d]&\Bp(\Kb)\ar@{->>}[d]_{\psi}&\Bp(\Kb)\times
\rH^1(B_{\etab},\F_p)\ar[rr]^-(0.5){(\bullet,\bullet)_B}&&\F_p\ar@{=}[d]& \rH^1(B_{\eta},\F_p)\ar[u]_{\alpha_2}\\
G^{a+}(\Kb)\ar@{^(->}[r]&G(\Kb)& G(\Kb)\times
G^\vee(\Kb)(-1)\ar[rr]^-(0.5){(\bullet,\bullet)_G}&&\F_p&G^\vee(\Kb)\ar@{^(->}[u]_{\gamma_2}.}\]
By functoriality, the homomorphism $\phi$ is adjoint to $\alpha_2$
and $\psi$ is adjoint to $\gamma_2$, \ie we have
$(\phi(x),f)_B=(x,\alpha_2(f))_A$ and
$(\psi(y),g)_B=(y,\gamma_2(g))_G$, for all $x\in\Ap(\Kb)$, $y\in
\Bp(\Kb)$, $f\in \rH^1(B_\etab,\F_p)$ and $g\in G^\vee(\Kb)$. By
\ref{comp-con}(ii), $\psi$ induces a surjective homomorphism
$\Bp^{a+}(\Kb)\ra G^{a+}(\Kb)$ for all $a\in \Q_{\geq 0}$. Now
statement   \ref{prop-BK}(iii) follows from \eqref{BK-fini-abel} and
Theorem \ref{Thm-BK-abel} applied to $B$.
\end{rem}

\section{Proof of Theorem \ref{main-thm}(i)}

\subsection{} For $r\in \Q_{>0}$, we denote by $\G_{a,r}$ the additive group scheme over $\Sb_r$ \eqref{notations}.
  Putting $t=1-{1/p}$,  we identify $\G_{a,t}$ with  an abelian $\fppf$-sheaf over $\Sb_1$
  by the canonical immersion $i: \Sb_t\ra
  \Sb_1$.
 Let $\rF:\G_{a,1}\ra \G_{a,1}$ be the Frobenius homomorphism, and $c$ a $p$-th root of $(-p)$.
 It is easy to check
  that the morphism $\rF-c:\G_{a,1}\ra \G_{a,1}$, whose cokernel is denoted by $\bP$, factorizes
through the canonical reduction morphism
 $\G_{a,1}\ra \G_{a,t}$, and we have an exact sequence of abelian
$\fppf$-sheaves on $\Sb_1$.
 \beq\label{syn-fppf}0\ra \F_p\ra \G_{a,t}\xra{\rF-c} \G_{a,1}\ra \bP\ra
 0.\eeq This exact sequence  gives (\ref{fc-suite}) after
restriction to
  the small \'etale topos $(\Xb_1)_{\et}$ for a   smooth  $S$-scheme $X$.

\begin{prop}\label{key-lemma}  Assume $p\geq  3$. Let $G$ be a truncated Barsotti-Tate group of level
$1$ over $S$, and $t=1-1/p$.  Then we have the equality
\[\dim_{\F_p}\rU^{e}G^\vee(\Kb)=\dim_{\F_p}\Ker \biggl(\Lie(\Gb^\vee_t)\xra{\rF-c}\Lie(\Gb^\vee_1)\biggr),\]
where $\rU^eG^\vee(\Kb)$ is the Bloch-Kato filtration
\eqref{filt-BK-fini}, and the morphism in the right hand side is
obtained
 by applying the functor $\Hom_{\Sb_1}(\Gb_1,\_)$ to the map
$\rF-c:\G_{a,t}\ra \G_{a,1}$.
\end{prop}

\subsection{} Before proving this proposition, we deduce first Theorem \ref{main-thm}(i).
 Let $G$ be a truncated Barsotti-Tate group of level $1$ and height $h$  over $S$ satisfying the
 assumptions of
\ref{main-thm}, $d$ be the dimension of $\Lie(G^\vee_s)$ over $k$, and
$d^*=h-d$.
 It follows from \ref{prop-BK}(iii) and \ref{key-lemma}  that
\[\dim_{\F_p}G^{\frac{e}{p-1}+}(\Kb)=h-\dim_{\F_p}\Ker \biggl(\Lie(\Gb^\vee_t)\xra{\rF-c}\Lie(\Gb^\vee_1)\biggr).\]
Since $\Lie(\Gb_1^\vee)$ is a free  $\cO_{\Sb_1}$-module of rank
$d^*$,  we obtain immediately Theorem \ref{main-thm}(i) by applying
Proposition \ref{prop-HW} to $\lam=c$ and $M=\Lie(\Gb_1^\vee)$.\\

The rest of this section is dedicated  to the proof of Proposition
\ref{key-lemma}.

\begin{lemma}\label{lem-inj} Let $G$ be  a Barsotti-Tate group of level 1 over $S$, $t=1-\frac{1}{p}$.
Then the morphism $\phi:\Ext^1(\Gb_1,\F_p)\ra
\Ext^1(\Gb_1,\G_{a,t})$ induced by the morphism $\F_p\ra \G_{a,t}$
in \eqref{syn-fppf} is injective.
\end{lemma}

\begin{proof}
  By (\cite{Il} Th\'eor\`eme 4.4(e)),
there exists a Barsotti-Tate group $H$ over $S$ such that we have an
exact sequence $0\ra G\ra H\xra{\times p}H\ra 0,$ which induces a
long exact sequence
\[\Ext^1(\Hb_1,\F_p)\xra{\times p} \Ext^1(\Hb_1,\F_p)\ra \Ext^1(\Gb_1,\F_p)\ra \Ext^2(\Hb_1,\F_p).\]
It is clear  that the multiplication by $p$ on $\Ext^1(\Hb_1,\F_p)$
is 0, and $\Ext^2(\Hb_1,\F_p)=0$ by Corollary \ref{ext2nul}; hence
the middle morphism in the exact sequence above  is an isomorphism.
Similarly, using $\Ext^2(\Hb_1,\G_{a,t})=0$ (\cite{BBM} 3.3.2), we
prove that the natural map $\Ext^1(\Hb_1,\G_{a,t})\ra
\Ext^1(\Gb_1,\G_{a,t})$ is an isomorphism. So we get a commutative
diagram
\[\xymatrix{\Ext^1(\Hb_1,\F_p)\ar[r]\ar[d]^{\phi_H}&\Ext^1(\Gb_1,\F_p)\ar[d]^{\phi}\\
\Ext^1(\Hb_1,\G_{a,t})\ar[r]&\Ext^1(\Gb_1,\G_{a,t}),}\] where the
horizontal maps are isomorphisms. Now it suffices to prove that
$\phi_H$ is injective.

Let $\bK$ be the $\fppf$-sheaf on $\Sb_1$ determined by the
following exact sequences:
\[0\ra \F_p\ra \G_{a,t}\ra \bK\ra 0; \quad \, 0\ra\bK\ra\G_{a,1}\ra \bP\ra 0.\]
Applying the functors $\Ext^q(\Hb_1,\_)$, we get
\begin{align*}\Hom(\Hb_1,\bK)\ra\Ext^1(\Hb_1,\F_p)\xra{\phi_H}
\Ext^1(\Hb_1,\G_{a,t});\\
0\ra \Hom(\Hb_1,\bK)\ra \Hom(\Hb_1,\G_{a,1})\ra\Hom(\Hb_1,\bP).
\end{align*}
Since $\Hom(\Hb_1,\G_{a,1})=0$ (\cite{BBM} 3.3.2), the injectivity
of $\phi_H$ follows immediately.
\end{proof}

\subsection{}\label{sub3.1}  Assume that $p\geq  3$. Let $G$ be a commutative finite and flat group
scheme killed by $p$ over $S$, and $0\ra G\ra A\ra B\ra 0$ be  a
resolution of $G$ by abelian schemes \eqref{imm-abel}. We denote
$\cP(B)=\Coker(\cO_{\Bb_1}\xra{\rF-c}\cO_{\Bb_1})$ \eqref{defn-cp},
and similarly for $\cP(A)$. According to  Remark \ref{rem-syn-BK},
we have an identification
\beq\label{3.7.2}\rH^0(B_{\esb},\rU^e\Psi^1_B)\simeq
\rH^0(B_{\esb},\cP(B))=\rH^0(\Bb_1,\cP(B))\eeq
 as submodules of $\rH^0(B_{\esb},\Psi^1_B)$;  in the last
 equality, we have identified the topos $(B_{\esb})_{\et}$ with $(\Bb_1)_{\et}$.
We denote
 \begin{align*}\Ker(B,\rF-c)&=\Ker\biggl(\rH^1(\Bb_1,\G_{a,t})\xra{\rF-c} \rH^1(\Bb_1,\G_{a,1})\biggr)\\
 &=\Ker\biggl(\Lie(\Bb_t^\vee)\xra{\rF-c}\Lie(\Bb_1^\vee)\biggr),\\
 \Ker(B,\delta_E)&=\Ker\biggl(\rH^0(\Bb_1,\cP(B))\xra{\delta_E}\rH^2(B_{\esb},\F_p)\biggr),\end{align*}
 where $\delta_E$ is the morphism defined in Proposition
 \ref{prop-fc}(2); we have also
 similar notations for $A$.
    Since the exact sequence (\ref{fc-coh-suite}) is
 functorial in $X$, we have a commutative diagram
\beq\label{fc-fini-abel}
\xymatrix{0\ar[r]&\rH^1(B_{\esb},\F_p)\ar[r]\ar[d]^{\beta_1}&\Ker(B,\rF-c)\ar[r]\ar[d]^{\beta_2}
&\Ker(B,\delta_E)\ar[r]\ar[d]^{\beta_3}&0\\
0\ar[r]&\rH^1(A_{\esb},\F_p)\ar[r]&\Ker(A,\rF-c)\ar[r]&\Ker(A,\delta_E)\ar[r]&0.}\eeq

\begin{lemma}\label{3.8} The assumptions are those of
\eqref{sub3.1}.

\emph{(i)} In  diagram \eqref{fc-fini-abel}, we have
\begin{align}\Ker\beta_1&=\Ker\biggl(\Ext^1(B_{\esb},\F_p)\ra
\Ext^1(A_{\esb},\F_p)\biggr)=(G_{\et})^\vee(\Kb)(-1);\label{3.8.1}\\
\Ker\beta_2&=\Ker\biggl(\Lie(\Gb_t^\vee)\xra{\rF-c}\Lie(\Gb_1^\vee)\biggr);\label{3.8.2}\\
\Ker\beta_3&=\rH^0(B_{\esb},\cP(B))\cap N(1)\subset
\rH^0(B_{\esb},\Psi^1_B),\label{3.8.3}
\end{align}where $N(1)$ is defined in  \eqref{suite-cycles}.

\emph{(ii)} We have the equality \beq\label{dim-equ}\dim_{\F_p}
\rU^eG^\vee(\Kb)=\dim_{\F_p}\Ker\beta_1+\dim_{\F_p}\Ker\beta_3.\eeq
In particular, we have \beq\label{dim-inequ}\dim_{\F_p}
\rU^eG^\vee(\Kb)\geq\dim_{\F_p}\Ker\biggl(\Lie(\Gb_t^\vee)\xra{\rF-c}\Lie(\Gb_1^\vee)\biggr),\eeq
Moreover, the equality holds in \eqref{dim-inequ}  if and only if
the morphism $\Coker\beta_1\ra \Coker\beta_2$ induced by diagram
\eqref{fc-fini-abel} is injective.
\end{lemma}

\begin{proof} (i) By Corollary
\ref{rem-im-ext}, we have  a canonical isomorphism
$\Ext^1(X_{\esb},\F_p)\simeq \rH^1(X_{\esb},\F_p)$ for $X=A$ or $B$.
Hence   formula \eqref{3.8.1} follows easily by applying the
functors $\Ext^q(\_,\F_p)$ to the exact sequence $0\ra G_{\esb}\ra
A_{\esb}\ra B_{\esb}\ra 0$. Applying the morphism of functors
$\rF-c:\Ext^i(\_,\G_{a,t})\ra \Ext^i(\_,\G_{a,1})$ to the exact
sequence $0\ra G\ra A\ra B\ra 0$, we obtain the following
commutative diagram
\[\xymatrix{0\ar[r] &\Lie(\Gb_t^\vee)\ar[r]\ar[d]^{\rF-c}&\Lie(\Bb_t^\vee)
\ar[r]\ar[d]^{\rF-c}&\Lie(\Ab_t^\vee)\ar[d]^{\rF-c}\\
0\ar[r]&\Lie(\Gb_1^\vee)\ar[r]&\Lie(\Bb_1^\vee)\ar[r]&\Lie(\Ab_1^\vee)},
\]
where we have used  \eqref{Lie-alg-fini} and \eqref{Lie-alg-abel}.
Formula \eqref{3.8.2} follows immediately from this diagram. For
 \eqref{3.8.3},  using  (\ref{3.7.2}),
we have a commutative diagram
\[\xymatrix{\Ker(B,\delta_E)\ar[r]^{(1)}\ar[d]_{\beta_3}&\rH^0(B_{\esb},\Psi^1_B)\ar[d]\\
\Ker(A,\delta_E)\ar[r]^{(2)}&\rH^0(A_{\esb},\Psi^1_A),}\] where the
maps $(1)$  and $(2)$ are injective. Hence we obtain
\begin{align*}\Ker \beta_3&=\Ker(B,\delta_E)\cap \Ker\biggl(\rH^0(B_{\esb},\Psi^1(B))\ra
\rH^0(A_{\esb},\Psi^1(A))\biggr)\\
&=\Ker(B,\delta_E)\cap N(1).
\end{align*}The morphisms $\delta_E$ and $d_2^{1,0}$ are compatible in
the sense of Proposition \ref{prop-fc}, and Proposition
\ref{prop-BK} implies that $\Ker(B,\delta_E)\cap
N(1)=\rH^0(B_{\esb},\cP(B))\cap N(1)$, which proves \eqref{3.8.3}.

(ii) By the isomorphism \eqref{4.9.1} and the surjectivity of the
morphism $u$ in \eqref{suite-cycles}, we have
\begin{align}\label{3.7.1}&\dim_{\F_p}\rU^eG^\vee(\Kb)=\dim_{\F_p}\biggl(G^\vee(\Kb)(-1)\cap
\rU^e\rH^1(B_{\etab},\Z/p\Z)\biggr)\nonumber\\
&=\dim_{\F_p}\biggl((G_{\et})^\vee(\Kb)(-1)\biggr)+\dim_{\F_p}\biggl(N\cap
\rH^0(B_{\esb},\rU^e\Psi^1(B))(-1)\biggr).\end{align} Then the
equality \eqref{dim-equ} follows from (i) of this lemma and
 \eqref{3.7.2}. The rest part of (ii)
follows immediately from diagram \eqref{fc-fini-abel}.
\end{proof}

\subsection{\textit{Proof of Proposition \ref{key-lemma}}}  We
choose a resolution $0\ra G\ra A\ra B\ra 0$ of $G$ by abelian
schemes \eqref{imm-abel}.
 By Lemma \ref{3.8}, we have to prove that if  $G$ is a truncated
 Barsotti-Tate group of level 1 over $S$, the morphism $\phi_{12}:\Coker\beta_1\ra
\Coker\beta_2$ induced by diagram \eqref{fc-fini-abel} is
injective.

By \ref{rem-im-ext}, we have
$\rH^1_{\et}(X_{\esb},\F_p)=\rH^1_{\et}(\Xb_1,\F_p)=\Ext^1(\Xb_1,\F_p)$
for $X=A$ or $B$. Thus the morphism $\beta_1$ is canonically
identified to the morphism $\Ext^1(\overline{B}_1,\F_p)\ra
\Ext^1(\Ab_1,\F_p)$ induced by the map $A\ra B$.
 Applying the functors $\Ext^i(\_,\F_p)$ to $0\ra \Gb_1\ra\Ab_1\ra \overline{B}_1\ra
 0$, we obtain a long exact sequence
 \[\Ext^1(\overline{B}_1,\F_p)\ra \Ext^1(\Ab_1,\F_p)\ra \Ext^1(\Gb_1,\F_p)\ra \Ext^2(\overline{B}_1,\F_p).\]
Since $\Ext^2(\overline{B}_1,\F_p)=0$ by  \ref{ext2nul}, we have
$\Coker\beta_1= \Ext^1(\Gb_1,\F_p)$. The commutative diagram
\[\xymatrix{0\ar[r]&\Ker(B,\rF-c)\ar[r]\ar[d]_{\beta_2}&\Lie(B_t^{\vee})\ar[r]^{\rF-c}\ar[d]_{\gamma}&\Lie(B_1^{\vee})\ar[d]\\
0\ar[r]&\Ker(A,\rF-c)\ar[r]&\Lie(A_t^\vee)\ar[r]^{\rF-c}&\Lie(A_1^\vee),}\]
induces a canonical morphism
$\psi:\Coker\beta_2\ra\Coker\gamma=\Ext^1(\Gb_1,\G_{a,t})$. Let
$\phi: \Ext^1(\Gb_1,\F_p)\ra \Ext^1(\Gb_1,\G_{a,t})$ be the morphism
induced by the map $\F_p\ra \G_{a,t}$ in (\ref{syn-fppf}). Then  we have the following commutative diagram 
\[\xymatrix{\Ext^1(\Gb_1,\F_p)=\Coker\beta_1\ar[r]^{\quad\quad\quad\phi_{12}}\ar[dr]_{\phi}&\Coker\beta_2\ar[d]^{\psi}\\
&\Ext^1(\Gb_1,\G_{a,t}).}\]  Now Lemma \ref{lem-inj} implies that
$\phi$ is injective, hence so is $\phi_{12}$. This completes the
proof.

\section{The canonical filtration in terms of  congruence groups}

\subsection{}\label{defn-cong} Recall the  following  definitions in \cite{SOS}.
 For any $\lam \in \cO_{\Kb}$, let $\cG^{(\lambda)}$  be the
group scheme $\Spec (\cO_{\Kb}[T,\frac{1}{1+\lambda  T}])$ with the
comultiplication given by $T\mapsto T\otimes 1+ 1\otimes T +\lambda
T\otimes T,$ the counit by $T=0$ and the coinverse by $T\mapsto
-\frac{T}{1+\lambda T}$. If $v(\lambda )\leq e/(p-1)$, we put
$$P_{\lambda}(T)=\frac{(1+\lambda T)^p-1}{\lambda^p}\in \cO_{\Kb}[T]$$
and let $\phi_{\lambda }: \cG^{(\lambda)}\ra \cG^{(\lam^p)}$ be the
morphism of $\cO_{\Kb}$-group schemes defined on the level of Hopf algebras by
$T\mapsto P_{\lam}(T)$. We denote by $G_\lam$ the kernel of
$\phi_{\lam}$, so we have $G_\lam=\Spec
\bigl(\cO_{\Kb}[T]/P_\lam(T)\bigr).$ We call it, following Raynaud, the
\emph{congruence group of level $\lam$}. It is a finite flat group
scheme over $\Sb=\Spec(\cO_{\Kb})$ of rank $p$.

\subsection{}\label{lem-theta-lam} For all $\lam\in \cO_{\Kb}$ with
$v(\lam)\leq e/(p-1)$, let $\theta_{\lam}: G_\lam\ra \mu_p=\Spec(\cO_{\Kb}[X]/(X^p-1))$ be the
homomorphism given on the level of Hopf algebras by $X\mapsto 1+\lam
T$. Then $\theta_\lam\otimes\Kb$ is an isomorphism, and if $v(\lam)=0$,
$\theta_\lam$ is an isomorphism.
 For all $\lam,\gamma\in \cO_{\Kb}$ with $v(\gamma)\leq v(\lam)\leq e/(p-1)$,
let  $\theta_{\lam, \gamma}:G_\lam\ra G_\gamma$ be the map defined by the
homomorphism of Hopf algebras $T\mapsto (\lam/\gamma)T$. We have
$\theta_\lam=\theta_\gamma\circ \theta_{\lam,\gamma}$.

\subsection{}  Let $\lam\in\OKb$ with
$v(\lam)\leq e/(p-1)$, $A$ be an abelian scheme over $S$.  We define
\beq\label{theta-abel}\theta_\lam(A):\Ext^1_{\Sb}(A,G_{\lam})\ra
 \Ext^1_{\Sb}(A,\mu_p)\eeq
to be the homomorphism induced by the canonical morphism
$\theta_\lam:G_\lam\ra \mu_p$, where, by abuse of notations, $A$
denotes also the inverse image  of $A$ over
$\Sb$, and $\Ext^1_{\Sb}$ means the extension in the category of
abelian $\fppf$-sheaves over $\Sb$.
Similarly, let $G$ be a commutative finite and flat group scheme
killed by $p$ over $S$; we define
\beq\label{theta-finite}\theta_{\lam}(G):\Hom_{\Sb}(G,G_{\lam})\ra
\Hom_{\Sb}(G,\mu_{p})=G^\vee(\Kb)\eeq
 to be the homomorphism induced by
 $\theta_{\lam}$. If $G=\Ap$, where $A$ is an $S$-abelian scheme, the natural exact sequence
 $0\ra \Ap\ra A\xra{\times p}A\ra 0$ induces
 a commutative diagram
 \beq\label{theta-abel-finite}\xymatrix{\Hom_{\Sb}(\Ap,G_\lam)\ar[r]\ar[d]_{\theta_{\lam}(\Ap)}
 &\Ext^1_{\Sb}(A,G_{\lam})\ar[d]^{\theta_\lam(A)}\\
 \Ap^\vee(\Kb)\ar[r]&\Ext^1_{\Sb}(A,\mu_p),}\eeq where horizontal maps are
 isomorphisms \eqref{ext-abel-global}.
 Hence,  $\theta_\lam(A)$ is canonically
identified to $\theta_\lam(\Ap)$.

\begin{lemma}\label{cong-indep}
Let $\lam,\gamma\in \OKb$ with $v(\gamma)\leq v(\lam)\leq e/(p-1)$, $G$ be
a commutative finite and flat group scheme killed by $p$ over $S$.

\emph{(i)}  $\theta_\lam(G)$  is injective.

\emph{(ii)} The image of $\theta_\lam(G)$  is contained in that of
$\theta_{\gamma}(G)$.

\emph{(iii)}  The image of $\theta_\lam(G)$   depends only on
$v(\lam)$, and it
 is invariant under the action of the Galois group $\Gal(\Kb/K)$.
\end{lemma}

\begin{proof}  We have a commutative diagram
\beq\label{6.4.1}\xymatrix{\Hom_{\Sb}(G_\lam^\vee,G^\vee )\ar[r]^{\theta_\lam(G)}\ar[d]_{(1)}&\Hom_{\Sb}(\Z/p\Z,G^\vee)\ar@^{=}[d]\\
\Hom_{\etab}(G_\lam^\vee,G^\vee)\ar[r]^{(2)}&\Hom_{\etab}(\Z/p\Z,G^\vee),}\eeq
where the horizontal maps are induced by
$\theta_{\lam}^\vee:\Z/p\Z\ra G_\lam^\vee$, and the vertical maps
 are induced by the  base change  $\etab\ra \Sb$.  Since
$\theta_{\lam}$ is an isomorphism over the generic point
(\ref{lem-theta-lam}), the map $(2)$ is an isomorphism.  Hence
statement (i) follows from the fact that $(1)$ is injective by the
flatness of $G$ and $G_\lam$.

Statement (ii) follows easily from  the existence of the morphism
$\theta_{\lam, \gamma}:G_\lam\ra G_\gamma$ with
$\theta_\lam=\theta_\gamma\circ\theta_{\lam,\gamma}$. The first part of (iii)
follows immediately  from (ii).    Any  $\sigma\in \Gal(\Kb/K)$
sends the image of  $\theta_\lam(G)$  isomorphically to the image of
$\theta_{\sigma(\lam)}(G)$,  which coincides with  the former by the
first assertion of (iii).
\end{proof}

\subsection{Filtration by congruence groups}\label{filt-cong-fini}
Let $a$ be a rational number with $0\leq a \leq e/(p-1)$, and $G$ be
a commutative finite and flat group scheme over $S$ killed by $p$.
We choose  $\lam\in \OKb$ with $v(\lam)=a$,  and denote by
$G^\vee(\Kb)^{[a]}$ the image of  $\theta_\lam(G)$. By Lemma
\ref{cong-indep}, $G^\vee(\Kb)^{[a]}$  depends only on $a$, and not
on the choice of $\lam$. Then $\bigl(G^\vee(\Kb)^{[a]}, a\in \Q\cap
[0, e/(p-1)])$ is an exhaustive decreasing filtration of
$G^\vee(\Kb)$  by $\Gal(\Kb/K)$-groups.

\subsection{ }\label{compare-AG} Let $\lam\in \OKb$  with $0\leq v(\lam)\leq e/(p-1)$,
$f:A\ra S$ be an abelian scheme, and $\fb:\Ab\ra\Sb$  its base
change by $\Sb\ra S$. In (\cite{AG} \S 6), Andreatta and Gasbarri
consider the homomorphism $\theta'_\lam(A):
\rH^1_{\fppf}(\Ab,G_\lam)\ra\rH^1_{\fppf}(\Ab,\mu_p)$ induced by
$\theta_\lam$, where by abuse of notation, $G_\lam$ denotes also the
$\fppf$-sheaf $G_\lambda$ restricted to $\Ab$. We have a commutative
diagram \beq\label{diag-6.6.1}
\xymatrix{\Ext^1_{\Sb}(A,G_\lam)\ar[r]^{\varphi(G_\lam)}\ar[d]_{\theta_\lam(A)} &\rH^1_{\fppf}(\Ab,G_\lam)\ar[d]^{\theta'_\lam(A)}\\
\Ext^1_{\Sb}(A,\mu_p) \ar[r]^{\varphi(\mu_p)}
&\rH^1_{\fppf}(\Ab,\mu_p),}\eeq where the horizontal arrows are the
homomorphisms \eqref{morph-ext-tor}.

\begin{lemma}\label{theta-AG} \emph{(i)} The homomorphisms $\varphi(G_\lam)$ and
$\varphi(\mu_p)$ in \eqref{diag-6.6.1} are isomorphisms. In
particular, the homomorphism $\theta'_\lam(A)$ is canonically
isomorphic to $\theta_\lam(A)$ \eqref{theta-abel}.

\emph{(ii)}   The canonical morphism $\rH^1_{\fppf}(\Ab,\mu_p)\ra
\rH^1(A_{\etab},\mu_p)$ is an isomorphism. Let
$\rH^1(A_{\etab},\mu_p)^{[v(\lam)]}$ be the image of
$\theta'_{\lam}(A)$ composed with this isomorphism. Then via the
canonical isomorphism $\rH^1(A_{\etab},\mu_p)\simeq \Ap^\vee(\Kb)$,
the subgroup $\rH^1(A_{\etab},\mu_p)^{[v(\lam)]}$ is identified to
$\Ap^\vee(\Kb)^{[v(\lam)]}$.
\end{lemma}

\begin{proof}
(i) For $H=G_\lam$ or $\mu_p$, the ``local-global'' spectral
sequence induces an exact sequence
\begin{equation}\label{formula-1}0\ra
\rH^1_{\fppf}(\Sb,R^0_{\fppf}\overline{f}_*(H_{\Ab}))\ra
\rH^1_{\fppf}(\Ab,H_{\Ab})\xra{\psi(H)}
 \rH^0_{\fppf}(\Sb,R^1_{\fppf}\overline{f}_*(H_{\Ab})).\end{equation}  By Prop.
\ref{cor-im-ext} and \eqref{ext-abel-global}, we have isomorphisms
$$\rH^0_{\fppf} (\Sb,R^1_{\fppf}\overline{f}_*(H_{\Ab}))\simeq
\rH^0_{\fppf}(\Sb,\cExt^1_{\Sb}(A,H))\simeq \Ext^1_{\Sb}(A,H).$$
Therefore, we obtain a homomorphism
$\psi(H):\rH^1_{\fppf}(\Ab,H_{\Ab})\ra\Ext^1_{\Sb}(A,H)$. We check
that the composed map $\psi(H)\circ\varphi(H)$ is the identity
morphism on $\Ext^1_{\Sb}(A,H)$; in particular, $\psi(H)$ is
surjective. By Prop. \ref{cor-im-ext}, we have also
$R^0_{\fppf}\overline{f}_*(H_{\Ab})=H_{\Sb}$; on the other hand, it
follows from (\cite{AG} Lemma 6.2) that
$\rH^1_{\fppf}(\Sb,H_{\Sb})=0$. Hence  $\psi(H)$ is injective by the
exact sequence \eqref{formula-1}, and $\varphi(H)$ and $\psi(H)$ are
both isomorphisms.

(ii) We have a commutative diagram
\[\xymatrix{\Ext^1_{\Sb}(A,\mu_p)\ar[r]^{\varphi(\mu_p)}\ar[d]_{(1)}&\rH^1_{\fppf}(\Ab,\mu_p)\ar[d]^{(2)}\\
\Ext^1_{\etab}(A_{\etab},\mu_p)\ar[r]&\rH^1(A_{\etab},\mu_p),}\]
where the vertical maps are base changes to the generic fibres, and
the horizontal morphisms are \eqref{morph-ext-tor}, which are
isomorphisms in our case by \eqref{rem-im-ext} and  statement (i).
The morphism $(1)$ is easily checked to be an isomorphism using
\eqref{ext-abel-global}, hence so is the morphism $(2)$. The second
part of  statement (ii) is a consequence of (i).
\end{proof}

The following proposition, together  with Proposition \ref{prop-BK},
implies Theorem \ref{main-thm2}.

\begin{prop} Let $G$ be a commutative finite and flat group scheme
over $S$ killed by $p$. Then, for all rational numbers $0\leq a\leq
e/(p-1)$,  we have $G^\vee(\Kb)^{[a]}=\rU^{pa}G^\vee(\Kb)$, where
 $\rU^{\bullet}G^\vee(\Kb)$ is the Bloch-Kato filtration
\eqref{filt-BK-fini}.
\end{prop}

\begin{proof}Let $0\ra G\ra A\ra B\ra 0$ be
 a resolution of $G$ by abelian schemes \eqref{imm-abel}. We
have, for all $\lam \in\OKb$ with $0\leq v(\lam)\leq e/(p-1)$, a
commutative diagram
 \[\xymatrix{0\ar[r]&\Hom_{\Sb}(G,G_\lam)\ar[r]\ar@{^(->}[d]^{\theta_\lam(G)}&\Ext^1_{\Sb}(B,G_\lam)
\ar[r]\ar@{^(->}[d]^{\theta_\lam(B)}&\Ext^1_{\Sb}(A,G_\lam)\ar@{^(->}[d]^{\theta_\lam(A)}\\
0\ar[r]&G^\vee(\Kb)\ar[r]&{}_pB^\vee(\Kb)\ar[r]&\Ap^\vee(\Kb).}\]
Hence,  for all rational numbers $a$ satisfying  $0\leq a\leq e/(p-1)$, we
have by \ref{theta-AG}(ii)
\beq\label{3-filtra}G^\vee(\Kb)^{[a]}=G^\vee(\Kb)\cap
{}_pB^\vee(\Kb)^{[a]}=G^\vee(\Kb)\cap \rH^1(B_{\etab},\mu_p)^{[a]}.
\eeq

 According to  (\cite{AG}
 Theorem 6.8),  the filtration $(\rH^1(B_\etab,\mu_p)^{[a]},0\leq a\leq e/(p-1))$
 coincides with the filtration $(\rU^{pa}\rH^1(B_{\etab},\mu_p),{0\leq a\leq \frac{e}{p-1}})$ (\ref{BK-abel}).
 Hence by (\ref{3-filtra}) and (\ref{BK-fini-abel}), the two filtrations $\bigl(G^\vee(\Kb)^{[a]}\bigr)$ and
 $\bigl(\rU^{pa}G^\vee(\Kb)\bigr)$ on $G^\vee(\Kb)$ coincide. This
 completes the proof.

\end{proof}

\section{The lifting property of the canonical subgroup}

In this section, by abuse of
 notations, $\G_a$ will denote  the additive group both over
$S$ and over $\Sb$.
 For a rational number $r>0$, we denote
 by $\G_{a,r}$, $\cG^{(\lam)}_r$ and $G_{\lam,r}$ the base changes to
 $\Sb_r$ of the respective group schemes.

\subsection{}\label{rel-TO} Following \cite{TO}, for $a,c\in \cO_{\Kb}$ with $ac=p$, we denote by
$G_{a,c}$ the group scheme $\Spec \bigl(\cO_{\Kb}[y]/(y^p-ay)\bigr)$ over
$\cO_{\Kb}$ with comultiplication
$$y\mapsto y\otimes 1+1\otimes y+
\frac{cw_{p-1}}{1-p}\sum_{i=1}^{p-1}\frac{y^i}{w_i}\otimes
\frac{y^{p-i}}{w_{p-i}}$$ and the counit given by $y=0$, where $w_i\,(1\leq
i\leq p-1)$ are universal constants in $\cO_{\Kb}$ with $v(w_i)=0$ (see
\cite{TO} p.9). Tate and Oort proved that $(a,c)\mapsto G_{a,c}$
gives a bijection between equivalence classes of factorizations of
$p=ac$ in $\cO_{\Kb}$ and isomorphism classes of $\cO_{\Kb}$-group schemes of
order $p$, where two factorizations $p=a_1c_1$ and $p=a_2c_2$ are
called equivalent if there exists $u\in \cO_{\Kb}$ such that
$a_2=u^{p-1}a_1$ and $c_2=u^{1-p}c_1$.

Let $\lam\in \cO_{\Kb}$ with $0\leq v_p(\lam)\leq 1/(p-1)$. There exists a
factorization $p=a(\lam)c(\lam)$ such that $G_\lam\simeq
G_{a(\lam),c(\lam)}$. More explicitly, we may take
$c(\lam)=\frac{(\lam(1-p))^{p-1}}{w_{p-1}}$ and
$a(\lam)=\frac{p}{c(\lam)}$, and we notice that
$v_p(a(\lam))=1-(p-1)v_p(\lam)$ is well defined independently of the
factorization $p=a(\lam)c(\lam)$.

\begin{lemma}[\cite{AG} Lemma 8.2 and 8.10]\label{AG8.2} Let $0<r\leq 1$ be a
rational number and $\lam\in\OKb$ with $v_p(\lambda)\leq 1-1/p$
 and $v_p(\lam^{p-1})\geq r$.

\emph{(i)}  Let $\rho^{\lam}_{r}$ be the morphism of groups schemes
${\cG}^{(\lam)}_r=\Spec\bigl(\cO_{\Sb_t}[T]\bigr)\ra
{\G}_{a,r}=\Spec\bigl(\cO_{\Sb_r}[X]\bigr)$ defined on the level of
Hopf algebras by $X\mapsto
\sum_{i=1}^{p-1}(-\lam)^{i-1}\frac{T^i}{i}$. Then $\rho^{\lam}_r$ is
an isomorphism. Moreover, the following diagram is commutative
\[\xymatrix{\cG^{(\lam)}_r\ar[r]^{\phi_{\lam,r}}\ar[d]_{\rho_r^\lam}&\cG_{r}^{(\lam^p)}\ar[d]^{\rho_r^{\lam^p}}\\
\G_{a,r}\ar[r]^{\rF-a(\lam)}&\G_{a,r},}\] where $\rF$ is the
Frobenius homomorphism and $a(\lam)\in \OKb$ is introduced in
\eqref{rel-TO}.

\emph{(ii)} Let $\delta_{\lam,r}$ be the composed morphism
$G_{\lam,r}\xra{i}\cG_r^{(\lam)}\xra{\rho^{\lam}_r}\G_{a,r}.$ Then
$\delta_{\lam,r}$ generates $\Lie(G_{\lam,r}^\vee)\simeq
\Hom_{\Sb_r}(G_{\lam,r},\G_{a,r})$ as an $\cO_{\Sb_r}$-module.
\end{lemma}

\begin{lemma}[\cite{AG} Lemma 8.3]\label{AG8.3} Let $\lam\in \OKb$
with $\frac{1}{p(p-1)}\leq v_p(\lam)\leq\frac{1}{p-1}$, and  $r
=(p-1)v_p(\lam)$. Then the following diagram is commutative
\[\xymatrix{\cG_1^{(\lam)}\ar[d]_{\iota_{1,r}}\ar[rr]^{\phi_{\lam,1}}&&\cG^{(\lam^p)}_1\ar[d]^{\rho^{\lam^p}_r}\\
\cG_r^{(\lam)}\ar[r]^{\rho_r^{\lam}}&\G_{a,r}\ar[r]^{\rF-a(\lam)}&\G_{a,1},}\]
where $\iota_{1,r}$ is the reduction map.
\end{lemma}

\subsection{}\label{assump7.3} Let   $\lam \in \OKb$ with
$v_p(\lam)=\frac{1}{p}$,
 $t=1-1/p$, and $G$ be a commutative finite and flat group
scheme killed by $p$ over $S$. We  define $\Phi_G$ to be
\begin{equation}\label{defn-Phi}
\Phi_G: \Hom_{\Sb}(\Gb,G_\lam)\xra{\iota_t}
\Hom_{\Sb_t}(\Gb_t,G_{\lam,t})\xra{\delta}\Hom_{\Sb_t}(\Gb_t,\G_{a,t})=\Lie(\Gb^\vee_t),
\end{equation}
where $\iota_t$ is the canonical reduction map, and $\delta$ is the
morphism induced by the element $\delta_{\lam,t}\in
\Hom_{\Sb_t}(G_{\lam,t},\G_{a,t})$ (\ref{AG8.2}(ii)).

\begin{prop}\label{prop-cong-fc} Let   $\lam \in \OKb$ with
$v_p(\lam)=\frac{1}{p}$,
 $t=1-1/p$, and $G$ be a  a Barsotti-Tate
group of level 1 over $S$, satisfying the hypothesis of Theorem
$\ref{main-thm}$. Then we have an exact sequence
\beq\label{fc-cong}0\ra
\Hom_{\Sb}(\Gb,G_\lam)\xra{\Phi_G}\Hom_{\Sb_t}(\Gb_t,\G_{a,t})\xra{\rF-a(\lam)}\Hom_{\Sb_{1}}(\Gb_1,\G_{a,1}),\eeq
 where $a(\lam)$ is defined   in
 \ref{rel-TO}.
\end{prop}

\begin{proof}
 From Lemma \ref{AG8.3}, we deduce a commutative diagram
\begin{equation}\label{7.4.1}\xymatrix{\Hom_{\Sb_1}(\Gb_1,G_{\lam,1})\ar[r]\ar[d]_{\iota_{1,t}}&\Hom_{\Sb_1}
(\Gb_1,\cG^{(\lam)}_1)\ar[rr]^{\phi_{\lam,1}}\ar[d]_{\iota'_{1,t}}
&&\Hom_{\Sb_1}(\Gb_1,\cG^{(\lam^p)}_1)\ar[d]^{\rho_1^{\lam^p}}\\
\Hom_{\Sb_t}(\Gb_t,G_{\lam,t})\ar[r]&\Hom_{\Sb_t}(\Gb_t,\cG^{(\lam)}_t)\ar[r]^{\rho_t^{\lam}}&\Hom_{\Sb_t}(\Gb_t,\G_{a,t})
\ar[r]^{\rF-a(\lam)}&\Hom_{\Sb_1}(\Gb_1,\G_{a,1}),}\end{equation}
where the upper row is exact and $\iota_{1,t}$ and $\iota'_{1,t}$
are reduction maps. Therefore, the composition of $\Phi_G$ with the
morphism $\rF-a(\lam)$ in (\ref{fc-cong}) factorizes through the
upper row of \eqref{7.4.1}, and equals thus 0. Let $L$ be the kernel
of the map $\rF-a(\lam)$ in (\ref{fc-cong}). Then $\Phi_G$ induces a
map $\Phi':\Hom_{\Sb}(\Gb,G_{\lam})\ra L$. We have to prove that
$\Phi'$ is an isomorphism.

Let $d^*$ be the rank of
$\Lie(\Gb_1^\vee)=\Hom_{\Sb_1}(\Gb_1,\G_{a,1})$ over $\cO_{\Sb_1}$,
and  recall that $v_p(a(\lam))=1/p$. Since $G$ satisfies the
assumptions of Theorem \ref{main-thm}, applying Prop. \ref{prop-HW}
to $\Lie(\Gb_1^\vee)$ and the operator $\rF-a(\lam)$, we see that
the group $L$ is an $\F_p$-vector space of dimension $d^*$. On the
other hand, $\Hom_{\Sb}(\Gb,G_{\lam})$ is identified with
$G^{\frac{e}{p-1}+}(\Kb)^{\perp}$ by Theorem \ref{main-thm2}. Thus
 it is also an $\F_p$-vector space of dimension $d^*$ by
Theorem \ref{main-thm}(i). Therefore, to finish the proof, it
suffices to prove that $\Phi'$ is surjective.

By Lemma \ref{AG8.2}(i), we have the following commutative diagram
\[\xymatrix{0\ar[r]&L\ar[r]\ar[d]_\alpha &\Hom_{\Sb_t}(\Gb_{t},\G_{a,t})\ar[r]^{\rF-a(\lam)}\ar@^{=}[d]
&\Hom_{\Sb_1}(\Gb_1,\G_{a,1})\ar[d]^{\iota_{1,t}}\\
0\ar[r]&\Hom_{\Sb_t}(\Gb_t,G_{\lam,t})\ar[r]&\Hom_{\Sb_t}(\Gb_t,\G_{a,t})\ar[r]^{\rF-a(\lam)}&\Hom_{\Sb_t}(\Gb_t,\G_{a,t}),}\]
where $\iota_{1,t}$ is the reduction map. The composed morphism
$\alpha\circ\Phi'$ is  the canonical reduction map, whose
injectivity will implies the injectivity of $\Phi'$. Thus the
following lemma will conclude the proof of the proposition.
\end{proof}

\begin{lemma} Assume that $p\geq 3$. Let $t=1-1/p$, $\lam\in \OKb$ with $v_p(\lam)=1/p$,
 and $G$ be a commutative finite  flat group scheme killed by $p$
over $S$. Then the  reduction map
\[\iota_t :\Hom_{\Sb}(\Gb,G_\lam)\ra \Hom_{\Sb_t}(\Gb_{t},G_{\lam,t})\]
is injective.
\end{lemma}
\begin{proof} We put $G\times_S\Sb=\Spec (A)$, where $A$ is a
Hopf algebra over $\OKb$ with the comultiplication $\Delta$. An
element $f\in \Hom_{\Sb}(\Gb,G_\lam)$ is determined by an element
$x\in A$ satisfying
\begin{align}\Delta(x)&=x\otimes 1+1\otimes x+\lam
x\otimes x\nonumber\\
P_\lam(x)&=\frac{(1+\lam x)^p-1}{\lam^p}=0.\label{7.3.1}\end{align}
Suppose that $\iota_t(f)=0$, which means $x\in \m_t A$. We want to
prove that in fact $x=0$. Let us write $x=\lam^ay$ where $a\geq p-1\geq 2$
is an integer, and $y\in A$. Substituting $x$ in (\ref{7.3.1}), we
obtain
\[(\lam^ay)^p+\sum_{i=1}^{p-1}\frac{1}{\lam^p}\binom{p}{i}\lam^{i(a+1)}y^i=0.\]
Since $v_p(\frac{1}{\lam^p}\binom{p}{i})=0$ for $1\leq i\leq p-1$ and $A$ is flat over $\cO_K$, 
we see easily that $y=\lam^{a+1}y_1$ for some $y_1\in A$. Continuing
this process, we find that $x\in \cap_{a\in\Q_{>0}}\m_a A=0$
\eqref{notations}.
\end{proof}

\begin{lemma}\label{lem-dim-Lie} Let $G$ be a Barsotti-Tate group of level $1$ and  height $h$ over $S$,
and $H$ be a flat closed subgroup scheme of $G$. We denote by $d$
the dimension of $\Lie(G_s)$ over $k$, and $d^*=h-d$. Then the
following conditions are equivalent:

\emph{(i)} The special fiber $H_s$ of $H$ coincides with the kernel
of the Frobenius of $G_s$.

\emph{(i')} The special fiber $H^\perp_s$ of $H^\perp=(G/H)^\vee$
coincides with the kernel of the Frobenius of $G^\vee_s$.

\emph{(ii)} $H$ has rank $p^d$ over $S$ and  $\dim_{k}\Lie
(H_{s})\geq d$.

 \emph{(ii')} $H^\perp$ has rank $p^{d^*}$ over $S$ and $\dim_{k}\Lie (H_{s}^\perp)\geq
d^*$.
\end{lemma}
\begin{proof} We have two exact sequences
\begin{align*}0&\ra H\ra G\ra G/H\ra 0\\0&\ra
H^{\perp}\ra G^\vee\ra H^\vee\ra 0.
\end{align*} Denote by $\ff_{G_s}$ (\resp by $\fV_{G_s}$) the Frobenius (\resp the Verschiebung) of $G_s$.
 Assume that (i) is satisfied, then we have $H^\vee_s=\Coker(\fV_{G_s^\vee})$ by duality. Since $G$ is a
Barsotti-Tate group of level 1,  $H^{\perp}_s$ coincides  with
$\im(\fV_{G_s^\vee})=\Ker (\ff_{G_s^\vee})$. Conversely, if
$H^{\perp}_s=\Ker(\ff_{G_s^\vee})$, we have also $H_s=\Ker
(\ff_{G_s})$. This proves the equivalence of (i) and (i').

If  (i) or (i') is satisfied, then (ii) and (ii') are also satisfied
(SGA $3_1$ $\mathrm{VII_{A}}$ 7.4). Assume (ii)  satisfied. Since
$\Ker(\ff_{H_s})$ has rank $p^{\dim_{k}\Lie(H_s)}$ (\emph{loc.
cit.}) and is contained in both $\Ker (\ff_{G_s})$ and $H_s$,
condition (ii) implies that these three groups have the same rank;
hence they coincide. This proves that  (ii) implies (i). The
equivalence of (i') and (ii') is proved in the same way.\end{proof}

 \subsection{\textit{Proof of Theorem $\ref{main-thm}(ii)$. }} By
 Lemma \ref{lem-dim-Lie},  the following lemma will complete the proof of
 \ref{main-thm}(ii).

\begin{lemma} Let $G$ be a  Barsotti-Tate group of level 1 and
  height $h$ over $S$,
satisfying the hypothesis of Theorem $\ref{main-thm}$, $d$ be the
dimension of $\Lie(G^\vee_s)$ over $k$, and $d^*=h-d$. Let $H$ be
the flat closed subgroup scheme $G^{\frac{e}{p-1}+}$, and
$H^\perp=(G/H)^\vee$. Then 
$H^\perp$  has  rank $p^{d^*}$ over $S$ and
$\dim_k\Lie(H^\perp_s)\geq d^*$.
\end{lemma}
 \begin{proof} Since $H$ has rank $p^d$ over $S$ by  \ref{main-thm}(i),
  $H^\perp$ has rank $p^{d^*}$ over $S$ and $\dim_{\F_p}(G/H)(\Kb)=d^*$. Let
$\lam\in\OKb$ with $v_p(\lam)=1/p$, and $t=1-1/p$. The canonical
projection $G\ra G/H$ induces an injective homomorphism
\begin{equation}\label{inj-G/H-G}\Hom_{\Sb}(\Gb/\Hb,G_\lam)\ra\Hom_{\Sb}(\Gb,G_\lam).\end{equation}
By Theorem \ref{main-thm2},
$H^\perp(\Kb)^{[\frac{e}{p}]}=\Hom_{\Sb}(\Gb/\Hb,G_{\lam})$ is
orthogonal to $(G/H)^{\frac{e}{p-1}+}(\Kb)$ under the perfect
pairing $(G/H)(\Kb)\times H^\perp(\Kb)\ra \mu_p(\Kb)$. As
$H=G^{\frac{e}{p-1}+}$, Prop. \ref{comp-con}(ii) implies that the
group scheme $(G/H)^{\frac{e}{p-1}+}$ is trivial. Hence we have
$$\dim_{\F_p}\Hom_{\Sb}(\Gb/\Hb,G_{\lam})=\dim_{\F_p}H^\perp(\Kb)=d^*=\dim_{\F_p}\Hom_{\Sb}(\Gb,G_\lam),$$
and  the canonical map \eqref{inj-G/H-G} is an isomorphism. By the
functoriality of $\Phi_G$ \eqref{defn-Phi}, we have a commutative
diagram
\begin{equation}\label{diag-Phis}\xymatrix{\Hom_{\Sb}(\Gb/\Hb,G_\lam)\ar@{=}[r]\ar[d]^{\Phi_{G/H}}&\Hom_{\Sb}(\Gb,G_\lam)\ar[d]^{\Phi_G}\\
\Lie(\overline{H}^\perp_{t})\ar[r]&\Lie(\Gb_t^\vee)},\end{equation}
where the lower row is an injective homomorphism of
$\cO_{\Sb_t}$-modules. Put $N_0=\Hom_{\Sb}(\Gb,G_\lam)$,
$M=\Lie(\Gb^\vee_1)$ and
$M_t=\Lie(\Gb_1^\vee)\otimes_{\cO_{\Sb_1}}\cO_{\Sb_t}=\Lie(\Gb^\vee_t)$.
By \ref{prop-cong-fc}(ii), $N_0$ is identified with the kernel of
$\rF-a(\lam):M_t\ra M$. Let $N$ be the $\cO_{\Kb}$-submodule of
$M_t$ generated by $N_0$. Applying \ref{prop-HW}(ii) to the morphism
$\rF-a(\lam)$, we get
$$\dim_{\kb}(N/\m_{\Kb}N)=\dim_{\F_p}N_0=d^*.$$  By
\eqref{diag-Phis},  $N$ is contained in $M'=\Lie(\Hb^\perp_t)\subset
M$. By applying  Lemma \ref{lem-ele} (ii) below to $N\subset M'$, we
obtain
\begin{equation}\label{inequality}d^*=\dim_{\kb}(N/\m_{\Kb}N)\leq
\dim_{\kb}(M'/\m_{\Kb}M').\end{equation}

Let $\omega_{\Hb^\perp_t}$ be the module of invariant differentials
of $\Hb^{\perp}_t$ over $\cO_{\Sb_t}$. Then we have $\omega_{H^\perp_{\esb}}=\omega_{\Hb^\perp_t}\otimes_{\cO_{\Sb_t}}\kb$ and 
$$M'=\Lie(\Hb^{\perp}_t)=\Hom_{\cO_{\Sb_t}}(\omega_{\Hb^\perp_t},\cO_{\Sb_t}).$$
Applying Lemma \ref{lem-ele} (i) to $\omega_{\Hb^\perp_t}$, we obtain
\begin{equation}\label{dim-equality-1}
\dim_{\kb}(M'/\m_{\Kb}M')=\dim_{\kb}\omega_{H^\perp_{\esb}}.
\end{equation}
From the relations $\Lie(H^\perp_s)\otimes_k \kb=\Lie(H^\perp_\esb)=\Hom_{\kb}(\omega_{H^\perp_\esb},\kb)$, we deduce 
\begin{equation}\label{dim-equality-2}
\dim_{\kb}\omega_{H^\perp_{\esb}}=\dim_k \Lie(H^\perp_s).
\end{equation}
The desired inequality $\dim_k \Lie(H^\perp_s)\geq d^*$ then follows from \eqref{inequality}, \eqref{dim-equality-1} and \eqref{dim-equality-2}.

\end{proof}

\begin{lemma}\label{lem-ele} Let $t$ be a positive rational number, $M$ be an  $\cO_{\Sb_t}$-module
of finite presentation.

\emph{(i)} Put $M^*=\Hom_{\cO_{\Sb_t}}(M,\cO_{\Sb_t})$. Then we have $\dim_{\kb}(M^*/\m_{\Kb}M^*)=\dim_{\kb}(M/\m_{\Kb}M)$.

\emph{(ii)} If  $N$ is a finitely presented
$\cO_{\Sb_t}$-submodule of $M$, then 
$\dim_{\kb}(N/\m_{\Kb}N)\leq \dim_{\kb}(M/\m_{\Kb}M)$.

\end{lemma}

\begin{proof}Since $M$ is of finite presentation, up to replacing $K$ by a finite extension, we may assume that there exists a positive integer $n$ and a finitely generated $\cO_K/\pi^n\cO_K$-module $M_0$, where $\pi$ is a uniformizer of $\cO_K$,  such that $\cO_{\Sb_t}=\cO_{\Kb}/\pi^n\cO_{\Kb}$ and  $M=M_0\otimes_{\cO_K}\cO_{\Kb}$. Note that there exist integers $0<a_1\leq \cdots a_r\leq n$  such that 
we have an exact sequence of $\cO_K$-modules
\begin{equation}\label{lem-ele-suite}0\ra \cO_{K}^r\xra{\varphi}\cO_{K}^r\ra M_0\ra 0,\end{equation}
where $\varphi$ is given by $(x_i)_{1\leq i\leq r}\mapsto (\pi^{a_i}x_i)_{1\leq i\leq r}$. 
 In order to prove (i), it suffices to verify that $\dim_k(M_0^*/\pi M_0^*)=\dim_k(M_0/\pi M_0)$, where $M_0^*=\Hom_{\cO_K}(M_0,\cO_K/\pi^n\cO_K)$.  Let 
\[(\cO_K/\pi^n\cO_K)^r\xra{\varphi_n}(\cO_K/\pi^n\cO_K)^r\ra M_0\ra 0\]
 be the reduction of \eqref{lem-ele-suite} modulo $\pi^n$.
Applying the functor $\Hom_{\cO_K}(\_, \cO_K/\pi^n\cO_K)$ to the above exact sequence,  we get
$$0\ra M_0\ra (\cO_{K}/\pi^n\cO_K)^r\xra{\varphi_n^*} (\cO_{K}/\pi^n\cO_K)^r$$ with $\varphi^*_n=\varphi_n$. Hence  $M_0^*$ is isomorphic to the submodule $\oplus_{i=1}^r(\pi^{n-a_i}\cO_{K}/\pi^n\cO_K)$ of $(\cO_{K}/\pi^n\cO_K)^r$, and we have $$\dim_k(M^*_0/\pi M_0^*)=r=\dim_k(M_0/\pi M_0).$$
 For statement (ii), by the same reasoning, we may assume that there exists a finite $\cO_{K}$-submodule $N_0$ of $M_0$ such that $N=N_0\otimes_{\cO_K}\cO_{\Kb}$.We need to prove that $\dim_{k}(N_0/\pi N_0)\leq \dim_k(M_0/\pi M_0)$. 
Let ${}_\pi M_0$ be the kernel of $M_0$ of the multiplication by $\pi$. We have an exact sequence of Artinian modules
\[0\ra {}_\pi M_0\ra M_0\xra{\times \pi}M_0\ra M_0/\pi M_0\ra 0.\]
By the additivity of length of Artinian modules, we obtain $\dim_k ({}_\pi M_0)=\dim_k(M_0/\pi M_0)$. Similarly, we have $\dim_k({}_\pi N_0)=\dim_k(N_0/\pi N_0)$. The assertion follows from the fact that ${}_\pi N_0$  is a submodule of $ {}_\pi M_0$.
\end{proof}

\begin{rem} If we could prove the exact sequence (\ref{fc-cong})
without knowing \emph{a priori} the rank of $\Hom_{\Sb}(\Gb,G_\lam)$
for $v_p(\lam)=1/p$, then we would get another proof of the
existence of the canonical subgroup of $G$. Since then, by
Proposition \ref{prop-HW} and (\ref{fc-cong}),
$\Hom_{\Sb}(\Gb,G_\lam)$ has $\F_p$-rank $d^*$ under the assumptions
of \ref{main-thm}. Then we identify it to be a subgroup of
$G^\vee(\Kb)$ by $\theta_\lam(G)$ (\ref{theta-finite}), and define
$H$ to be the subgroup scheme of $G$ determined by
$H(\Kb)^{\perp}=\Hom_{\Sb}(G,G_\lam).$ The arguments in this section
imply that $H$ is the canonical subgroup of $G$. For abelian
schemes, this approach is due to Andreatta-Gasbarri \cite{AG}.
\end{rem}

\

\begin{thebibliography}{99}
\bibitem{AM} \textsc{A. Abbes} and \textsc {A. Mokrane},
Sous-groupes canoniques et cycles \'evanescents $p$-adiques pour les
vari\'et\'es ab\'eliennes, \emph{Publ. Math. Inst. Hautes \'Etud.
Sci.} \textbf{99} (2004), 117-162.

\bibitem{AS} \textsc{A. Abbes} and \textsc{ T. Saito}, Ramification
of local fields with imperfect residue fields, \emph{Am. J. Math.}
\textbf{124} (2002), 879-920.

\bibitem{AG} \textsc{F. Andreatta} and \textsc{C. Gasbarri}, The
canonical subgroup for families of abelian varieties, \emph{Compos. Math.} \textbf{143} (2007), 566-602.


\bibitem{Ar} \textsc{M. Artin}, Algebraization of formal moduli, \emph{Global analysis},
Univ. Tokyo Press  (1969), 21-71.

\bibitem{BBM} \textsc{F. Berthelot, L. Breen} and \textsc{W.
Messing}, \emph{Th\'eorie de Dieudonn\'e Cristalline II}, LNM
\textbf{930}, Springer-Verlag, (1982).

\bibitem{BK} \textsc{S. Bloch} and \textsc{K. Kato}, $p$-adic
\'etale cohomology, \emph{Publ. Math. Inst. Hautes \'Etud. Sci.}
\textbf{63} (1986), 107-152.



\bibitem{Co} \textsc{B. Conrad}, Higher-level canonical subgroups in
abelian varieties, Preprint (2005).

\bibitem{Dem} \textsc{M. Demazure},  Bidualit\'e des sch\'emas ab\'eliens, in \emph{S\'eminaire de g\'eom\'etrie alg\'ebrique}, Orsay, (1967-1968).

\bibitem{dwork} \textsc{B. Dwork}, $p$-adic cycles,
\emph{Publ. Math. Inst. Hautes \'Etud. Sci.} {\bf 37} (1969),
27-115.

\bibitem{FC} \textsc{G. Faltings} and \textsc{L. Chai}, \emph{Degeneration of abelian varieties},
Springer-Verlag (1990).

\bibitem{Gr} \textsc{A. Grothendieck}, Le groupe de Brauer III, in
\emph{Dix expos\'es sur la cohomologie des sch\'emas}, North-Holland
(1968).

\bibitem{Il} \textsc{L. Illusie}, D\'eformations de groupes de
Barsotti-Tate (d'apr\`es A. Grothendieck), \emph{Ast\'erisque}
\textbf{127} (1985), 151-198.

\bibitem{Ka} \textsc{K. Kato}, On $p$-adic vanishing cycles (Applications of ideas of
Fontaine-Messing),  \emph{Adv. Stud. Pure Math.}, \textbf{10}
(1987), 207-251.

\bibitem{Kz} \textsc{N. Katz}, $p$-adic properties of modular
schemes and modular forms, in \emph{Modular functions of one
variable III}, LNM \textbf{350}, Springer-Verlag, (1973).


\bibitem{KL} \textsc{M. Kisin} and \textsc{K. F. Lai},
Overconvergent Hilbert modular forms, \emph{Amer. J. of Math.}
\textbf{127} (2005), 735-783.

\bibitem{lubin}\textsc{J. Lubin}, Finite subgroups and isogenies of
one-parameter formal groups, \emph{Ann. of Math.} \textbf{85}, 2nd
series (1967), 296-302.

\bibitem{MM}\textsc{B. Mazur} and \textsc{W. Messing}, \emph{Universal
extensions and one dimensional crystalline cohomology}, LNM
\textbf{370}, Springer-Verlag, (1974).

\bibitem{Mu}\textsc{D. Mumford}, \emph{Geometric invariant theory},
Springer-Verlag, (1965).



\bibitem{Ray}\textsc{M. Raynaud}, Sp\'ecilisation du foncteur de
Picard, \emph{Publ. Math. Inst. Hautes \'Etud. Sci.} \textbf{38}
(1970),  27-76.



\bibitem{SOS} \textsc{T. Sekiguchi, F. Oort} and \textsc{N. Suwa},
On the deformation of Artin-Schreier to Kummer, \emph{Ann. Sci. de
l'\'E.N.S.} $4^e$ s\'erie, tome 22, No.3 (1989), 345-375.

\bibitem{TO} \textsc{J. Tate} and \textsc{F. Oort}, Group schemes of prime order,
 \emph{Ann. Sci. de
l'\'E.N.S.} $4^e$ s\'erie, tome 3, No. 1 (1970), 1-21.

\end{thebibliography}
\end{document}